\documentclass{amsart}

\usepackage{amsmath,amssymb}
\usepackage{bm}
\usepackage{graphicx}
\usepackage{ascmac}
 \usepackage{amsthm}
 
\theoremstyle{plain}
\newtheorem{thm}{Theorem}[section]

\newtheorem{prop}[thm]{Proposition}
\newtheorem{lmm}[thm]{Lemma}

\theoremstyle{definition}
\newtheorem{dfn}[thm]{Definition}
\newtheorem{rmk}[thm]{Remark}

\def\Spec{\mathop{\mathrm{Spec}}\nolimits}
\def\Proj{\mathop{\mathrm{Proj}}\nolimits}
\def\an{\mathop{\mathrm{an}}\nolimits}
\def\can{\mathop{\mathrm{can}}\nolimits}
\def\ord{\mathop{\mathrm{ord}}\nolimits}
\def\ddc{\mathop{\mathrm{dd^c}}\nolimits}
\def\MA{\mathop{\mathrm{MA}}\nolimits}
\def\hyb{\mathop{\mathrm{hyb}}\nolimits}
\DeclareMathOperator{\Lyap}{Lyap}
\DeclareMathOperator{\pro}{pr}
\DeclareMathOperator{\Diff}{D}

\newcommand{\bP}{\mathbb{P}}
\newcommand{\pr}{\mathbb{P}}

\newcommand{\bC}{\mathbb{C}}
\newcommand{\C}{\mathbb{C}}
\newcommand{\D}{\mathbb{D}}
\newcommand{\bD}{\mathbb{D}}
\newcommand{\N}{\mathbb{N}}

\newcommand{\bR}{\mathbb{R}}
\newcommand{\R}{\mathbb{R}}
\newcommand{\A}{\mathbb{A}}
\newcommand{\bA}{\mathbb{A}}
\newcommand{\Sh}{\mathcal{O}}
\newcommand{\F}{\mathcal{F}}

\usepackage[truedimen,top=25truemm,bottom=25truemm,left=25truemm,right=25truemm]{geometry}

\title{Non-archimedean and hybrid dynamics of H\'enon mappings}
\author{Reimi Irokawa}
\address{NTT Institute for Fundamental Mathematics, NTT Communication Science Laboratories, Nippon Telegraph and Telephone Corporation, 3-9-11 midori-cho, Musashino-shi, Tokyo 180-8686, Japan}
\email{reimi.irokawa@ntt.com}

\date{2023/4/6}

\begin{document}

\maketitle

\begin{abstract}
     For studying the meromorphic degeneration of complex dynamics, the theory of hybrid spaces, introduced by Boucksom, Favre and Jonsson, is known to be a strong tool. In this paper, we apply this theory to the dynamics of H\'enon maps. For a family of H\'enon maps $\{H_t\}_{t\in\D^*}$ that is parametrized by a unit punctured disk and meromorphically degenerates at the origin, we show that as $t\to 0$, the family of the invariant measures $\{\mu_t\}$ ``weakly converges'' to a measure on the Berkovich affine plane associated to the non-archimedean H\'enon map determined by the family $\{H_t\}_t$. We also calculate the limit of their Lyapunov exponents.
\end{abstract}

\section{Introduction}\label{sec:intro}
The dynamics of H\'enon maps were first studied by H\'enon, and later by other researchers, 
because quite rich and complicated dynamics can be observed despite the simple description. A generalized complex H\'enon map of degree $d\geq2$, which is the main target of this paper,  is defined as follows:
\begin{equation}
H(x,y)=(p(x)-ay,x),
\end{equation}
where $a\in\C\setminus\{0\}$, and $p(x)$ is a monic polynomial of degree $d$ with coefficients in $\C$. Bedford and Smillie studied the potential-theoretic properties of such a map \cite{BS1} and \cite{BS2}. More precisely, for a generalized H\'enon map $H(x,y)$ of degree $d$, they defined the following Green functions
\begin{align*}
    G^+(x,y)&=\lim_{n\to\infty}\frac{1}{d^n}\log^+\|H^{\circ n}(x,y)\|, \text{ and}\\
    G^-(x,y)&=\lim_{n\to\infty}\frac{1}{d^n}\log^+\|H^{\circ-n}(x,y)\|,  
\end{align*}
where $\|\cdot\|$ is any norm on $\C^2$. These are well-defined continuous functions, and the limits are locally uniform. By taking their Laplacians, we obtain the positive closed $(1,1)$-currents $T^+=\ddc G^+/2\pi$ and $T^-=\ddc G^-/2\pi$. These two currents satisfy $H^*(T^+)=dT^+$ and $H^*(T^-)=d^{-1}T^-$ respectively by definition. Then, as they both have continuous potentials, their wedge product yields a probability measure $\mu=T^+\wedge T^-$, which is invariant under the action of $H^*$. It is also known from \cite{BS1} that $\mu=\ddc(\max\{G^+,G^-\})^{\wedge 2}/4\pi^2$.

In this paper, we investigate degenerating families of these dynamics. Let $\mathcal{O}(\mathbb{D})[t^{-1}]$ be the set of meromorphic functions over the unit disk $\D$ with poles only at the origin $0$. Consider a family of generalized complex H\'enon maps $\{H_t\}_{t\in\D^*}$ that are parametrized by the punctured unit disc $\D^*=\{t\in\C\ |\  |t|<1, t\neq0\}$, i.e.,
\begin{equation}
H_t(x,y)=(p_t(x)-a(t)y,x),
\end{equation}
where $a(t)\in \mathcal{O}(\mathbb{D})[t^{-1}]$ without any zero on $\bD^*$, and $p_t(x)$ is a monic polynomial with coefficients in $ \mathcal{O}(\mathbb{D})[t^{-1}]$. For each $t\in\D^*$, $H_t$ defines a generalized H\'enon map. This paper's goal is to study the motion of the dynamics of such families as $t\to0$.

We apply the so-called theory of hybrid spaces for analyzing such families. Hybrid spaces were originally invented in algebraic geometry by Boucksom, Favre, and Jonsson \cite{BJ17}, \cite{BFJ15} and \cite{BFJ16} to study the degeneration of algebraic varieties. Later, Favre applied it the theory of hybrid spaces the degeneration of dynamical systems \cite{Fa}. Roughly speaking, the theory provides a way to treat the non-archimedean dynamics as the ``limit'' of degenerating families. 

Before introducing hybrid spaces, we study the non-archimedean dynamics. Note that the space $\mathcal{O}(\mathbb{D})[t^{-1}]$ can be naturally regarded as a subset of the field $\C((t))$. Then the family $\{H_t\}_{t\in\D^*}$ defined above naturally induces a polynomial automorphism of a scheme $\A^2_{\C((t))}=\Spec\C((t))[x,y]$ by regarding the $a_i(t)$ and $a(t)$ as elements of $\C((t))$. By taking the analytification of this family (in the sense of Berkovich \cite{Ber1}), we obtain the following automorphism:
\begin{align*}
    H:\A^{2,\an}_{\C((t))}\to\A^{2,\an}_{\C((t))}.
\end{align*}

In Section \ref{sec:non-archimedean-henon-maps}, we investigate the properties of this map, a non-archimedean H\'enon map. In particular, for such $H$, the Green functions $G^+$ and $G^-$ can be defined in just the same way as for the complex H\'enon map (see Proposition \ref{prop:thenonarchimedeangreenfunction} below). From the theory of non-archimedean Monge-Amp\`ere operator develloped by Boucksom, Favre and Jonsson in papers \cite{BFJ15} and \cite{BFJ16}, or the theory of forms and currents over Berkovich spaces by Chambert-Loir and Ducros in \cite{CLD}, we introduce a positive measure $\mu=\ddc(\max\{G^+,G^-\})^{\wedge 2}$. Note that there has been little research on the non-archimedean dynamics of H\'enon maps. In \cite{DP}, DeMark and Petsche studied the horseshoes of non-archimedean H\'enon maps, but without using Berkovich spaces. Our work is related but entails a different situation.

Next, we introduce the notions of dynamics over hybrid spaces. Roughly spealing, the hybrid affine plane $\A^2_{\hyb}$ is the family that results from the complex affine plane $\C^2$ parametrized by $\D^*$ ``degenerating" to the Berkovich affine plane $\A^{2,\an}_{\C((t))}$ over $\C((t))$ at the origin. Indeed, after fixing a real number $0<r<1$, let $\bar{\D}_r$ denote the complex closed disk $\{t\in\C|\ |t|\leq r\}$ of radius $r$ centered at the origin. The space $\A^2_{\hyb}$ has a proper surjection $\pi:\A^2_{\hyb}\to\bar{\D}_r$ for some $0<r<1$ such that there exists a family of homeomorphisms $\{\psi_t\}_{t\in\bar{\D}_r}$, where
\begin{align*}
    &\psi_t:\C^2\to\pi^{-1}\{t\}\text{ for }t\neq0,\\
    &\psi_0:\A^{2,\an}_{\C((t))}\to\pi^{-1}\{0\}.
\end{align*}
For a brief review of the hybrid space's construction, see Section \ref{sec:mongeamperemeasuresonthehybridprojectiveplane}. On the hybrid affine plane, the family $\{H_t\}$ naturally induces an automorphism of $\A^2_{\hyb}$ that preserves the fiber structure by $\pi$. On the fiber $\pi^{-1}\{t\}$ for $t\neq0$, the morphism is exactly $H_t$, and on $\pi^{-1}\{0\}$, it is the $H$ defined above, via the homeomorphisms $\psi_t$. By leveraging these notions, we can regard the non-archimedean dynamics by $H$ as the ``limit" of the family of dynamics $\{H_t\}_{t\in\D^*}$. Then, under these settings, we have the following theorem, which is our main result:
\begin{thm} \label{thm:intromain}
 Let $\{H_t\}_{\D^*}$ be an analytic family of H\'enon maps, and let $H$ be the non-archimedean H\'enon map induced from $\{H_t\}$ as above. Let $\mu_t=\ddc(\max\{G_t^+,G_t^-\})^{\wedge 2}/4\pi^2$ for $t\neq0$, where $G^{\pm}_t$ are the Green functions attached to $H_t$, and let $\mu=\ddc(\max\{G^+,G^-\})^{\wedge 2}$, where $G^{\pm}$ are the (non-archimedean) Green functions attached to $H$. Then, the following weak convergence holds:
 \begin{align*}
     \lim_{t\to0,t\neq0}(\psi_t)_*\mu_t=(\psi_0)_*\mu_0.
 \end{align*}
\end{thm}
In Section \ref{sec: Aoneparameterfamilyofcomplexhenonmaps}, we explain the key lemma for proving this theorem, and then we complete the proof in Section \ref{sec:maintheorem}.

By means of this theorem, we also study the limit of Lyapunov exponents. For a generalized complex H\'enon map $H(x,y)=(p(x)-ay,x)$, its Jacobi matrix $\Diff H$ is written as 
\begin{align} \label{eq:DH}
    \Diff H= \begin{pmatrix}
    p'(x) & a \\
    1 & 0
    \end{pmatrix} .
\end{align}
For a non-archimedean H\'enon map $H$, we can consider the same matrix $\Diff H$, too.

We have several kinds of Lyapunov exponents: the average by the probability measure $\mu$ of the first eigenvalue of $DH$ is called the first Lyapunov exponent, that of the second one is the second Lyapunov exponent, and the sum of them is called the total Lyapunov exponent respectively. For a degenerating family $\{H_t\}_{t\in\bD^*}$ that we consider above, we denote them by $\lambda_{1,t}$, $\lambda_{2,t}$ and $\Lyap(H_t)$ respectively.  Then we have the following theorem:

\begin{thm}\label{thm:convoflyapunocexponents}
    Let $\{H_t\}_{t\in\bD^*}$, $\lambda_{1.t}$, and $\Lyap(H_t)$ be as above. Then, for the total Lyapunov exponents, we have
    \begin{align*}
    \Lyap(H_t)=\frac{\log|t|^{-1}}{\log r}\log|a(t)|+o(\log|t|^{-1}),
\end{align*}
where $|a(t)|$ is the norm of $a(t)\in\bC((t))$ evaluated by the $t$-adic norm to base $r$. Also, for the first Lyapunov exponents, we have
    \begin{align*}
        \lambda_{1,t}=\frac{\log|t|^{-1}}{\log r}\Lambda+o(\log|t|^{-1}),
    \end{align*}
    where
    \begin{align*}
        \Lambda=\lim_{n\to\infty}\frac{1}{n}\int_{\bA^{2,\an}_{\bC((t))}}\log\|\Diff H^n\|d\mu.
    \end{align*}
\end{thm}

As noted above, the theory of hybrid spaces was first applied by Favre to study the degeneration of dynamics \cite{Fa}. He showed a result similar to Theorem \ref{thm:intromain} for degenerating families of endomorphisms over projective spaces $\pr^N$ for $N\geq1$. The basic argument of Theorem \ref{thm:intromain} is parallel to that result: the essential part is to show the ``uniformity'' of the family of functions $\{\max\{G_t^+,G_t^-\}\}_{t\in\D^*}$. In this paper, we apply Favre's result to an automorphism of $\A^2_k$, which is not projective, or even proper. Indeed, the H\'enon map cannot be extended as an endomorphism of any proper variety, which is a significant obstacle to constructing an argument similar to that in \cite{Fa}. Favre showed a weak convergence similar to that in Theorem \ref{thm:intromain} but under a more general setting, as given in Theorem 4.2 of \cite{Fa}, and he then applied it to endomorphisms of $\pr^N(\C)$. Favre's Theorem 4.2 actually requires the space to be projective. In contrast, the proof here uses  $(H_t,H_t^{-1}):\C^2\to\C^4$ instead of $H_t$, which can be extended to a morphism $\pr^2(\C)\to\pr^4(\C)$ and thus enable application of Favre's result as explained in Section \ref{sec:maintheorem}. This kind of argument is done by Lee \cite{LeeG} to study number-theoretic properties of the so-called algebraic stable pair of polynomials, which is a kind of generalization of H\'enon maps to morphisms of higher dimensional affine spaces.



\subsection*{Notation}
As noted above, denote the set of meromorphic functions on the unit disk $\D^*$ with poles only at $0$ by $\mathcal{O}(\mathbb{D})[t^{-1}]$. For the derivative $\ddc$, we put
\[
\mathop{\mathrm{d^c}}:=-i(\partial-\bar{\partial}),
\]
so that $\ddc=2i\partial\bar{\partial}$.

\subsection*{Acknowledgement} I would like to thank my advisor, Prof.\ Fumiharu Kato for his kind advice. This research was supported by a JSPS Grant-in-Aid for JSPS Fellows, 20J14309.

\section{Non-archimedean H\'enon maps}\label{sec:non-archimedean-henon-maps}
Let $k$ be an algebraically closed field with complete, non-archimedean valuation. 
 In this section, we consider an automorphisms of $\A_k^2$ with the following form:
 \begin{align*}
    H(x,y)=(p(x)-ay,x),
\end{align*}
where $a,b\in k$, and $p=x^d+a_1x^{d-1}+\cdots+a_d$ is a monic polynomial of degree $d$ greater than $1$ with coefficients in $k$. Note that the inverse is
\begin{align*}
    H^{\circ -1}(x,y)=\left(y, \frac{p(y)-x}{a}\right).
\end{align*}

Now, we take the Berkovich analytification of $\A^2_k$, $H$, and $H^{\circ -1}$. By abuse of notation, we also use $H$ and $H^{\circ -1}$ to denote their respective analytifications. We have two analytic morphisms:
\begin{align*}
    H&:\A^{2,\an}_k\to\A^{2,\an}_k, \text{ and} \\
    H^{\circ -1}&:\A^{2,\an}_k\to\A^{2,\an}_k. \\
\end{align*}

In this paper, we study the potential-theoretic properties of these maps. 

Recall that the underlying set of $\A^{2,\an}_k$ consists of the multiplicative seminorms on $k[x,y]$ whose restriction to $k$ coincides with the fixed norm. For any $f\in k[x,y]$ and $z\in\A^{2,\an}_k$, denote the value of $f$ at $z$ by $|f(z)|$. Note that if a point $z_0$ is a classical point i.e.\ $z_0=(x_0,y_0)\in k^2$, then $|f(z_0)|=|f(x_0,y_0)|$; that is, the corresponding seminorm is the norm of the evaluation map of $f$ at $(x_0,y_0)$. The topology of $\A^{2,\an}$ is the weakest one such that $z\mapsto |f(z)|$ is continuous for any $f\in k[x,y]$.

Later arguments will frequently use the following two continuous functions: $|x(z)|$ and $|y(z)|$. Let $|\pi_1(z)|:=|x(z)|$ and $|\pi_2(z)|=|y(z)|$, respectively.

Note that for $z_0=(x_0,y_0)\in \A^{2,\an}_k(k)$, $|\pi_1(z_0)|=|x_0|$ and $|\pi_2(z_0)|=|y_0|$, i.e., the functions $|\pi_1(z)|$ and $|\pi_2(z)|$ are the absolute values of the projections when restricted to rational points. Actually, outside $\A^{2,\an}_k(k)$, $|\pi_1(z)|$ and $|\pi_2(z)|$ can also be regarded as the norms of the natural projection at $z$. Indeed, let $\pi_1$ (resp. $\pi_2$) be the analytification of the natural projection $\A^2_k\to\A^1_k$, $(x,y)\mapsto x$ (resp.  $(x,y)\mapsto y$). Note that the projections comes from the ring homomorphism $x:k[x]\to k[x,y]$ (resp. $y:k[x]\to k[x,y]$). Then, on $\A^1_k=k[x]$ (resp. $\A^1_k=k[y]$), there is a regular function $x$ (resp. $y$), and we can consider the norm of $x$ (resp. $y$) evaluated at any point in $\A^{1,\an}_k$ just as in the two-dimensional case. For any $z\in\A^{2,\an}_k$, the value $|\pi_1(z)|$ (resp. $|\pi_2(z)|$) defined above actually coincides with $|x(\pi_1(z))|$ (resp. $|y(\pi_2(z))|$).

Next, define a real-valued function $\|\cdot\|$ on $\A^{2,\an}_k$ by extending the norm on $k^2$ as follows:
\begin{align*}
    \|z\|:=\max\{|\pi_1(z)|,|\pi_2(z)|\}.
\end{align*}

This section's goal is to prove the following theorem.
\begin{thm} \label{thm: non-archimedean-Green-function}
Define the functions $G^+$ and $G^-:\A^{2,\mathop{\mathrm{an}}}_k\to\R_{\geq0}$ as follows:
\begin{align*}
    G^+(z)&=\lim_{n\to\infty}\frac{1}{d^n}\log^+\|H^{\circ n}(z)\|, \text{ and}\\
    G^-(z)&=\lim_{n\to\infty}\frac{1}{d^n}\log^+\|H^{\circ-n}(z)\|,   
\end{align*}
where $\log^+r=\log\max\{1,r\}$ for $r\in\bR$. The functions $G^{\pm}$ are well-defined and continuous, and the above limits are both uniform on each compact set in $\A^{2, \an}_k$.
\end{thm}

Note that this statement is well-known from \cite{BS1} for the complex dynamical situation. The basic strategy here is parallel to the strategy there but requires addressing a few topological issues, as Berkovich spaces are not first-countable in general.

First, as in the complex case, several notations and lemmas are necessary. For any $R\in\bR$ sufficiently large that $R>\max\{|a_1|,\ldots,|a_d|,|a|,1\}$, we consider the following three sets:
\begin{align*}
    &V^+_R=\overline{\left\{(x,y)\in k^2\ \middle|\  |x|\geq|y|,\  |x|\geq R\right\}}, \\
    &V^-_R=\overline{\left\{(x,y)\in k^2\ \middle|\  |y|\geq|x|,\  |y|\geq R\right\}}, \\
    &W_R=\overline{\left\{(x,y)\in k^2\ \middle|\  |x|\leq R,\  |y|\leq R\right\}}, \\
\end{align*}
where, for a subset $E\subset\A^{2,\an}_k$, $\overline{E}$ denotes the topological closure in $\A^{2,\an}_k$. As $R$ can be made arbitrarily largely, the family $\{W_R\}_{R>\max\{|a_1|,\ldots,|a_d|,|a|,1\}}$ covers the whole $\A^{2,\an}_k$.

We also consider the following four sets:
\begin{align*}
    &K^+=\left\{z\in\A^{2, \an}\ \middle|\ \{\|H^{\circ n}(z)\|\}_n\text{ is bounded}\right\}, \\
    &K^+=\left\{z\in\A^{2, \an}\ \middle|\ \{\|H^{\circ -n}(z)\|\}_n\text{ is bounded}\right\}, \\
    &J^+=\partial K^+, \\
    &J^-=\partial K^-.
\end{align*}

There are several situations in which Theorem \ref{thm: non-archimedean-Green-function} is trivial, and we address those first.
\begin{lmm}\label{lmm:NAtrivcase}
For any $R>\max\{|a_1|,\ldots,|a_d|,|a|,1\}$ and $z\in V^+_R$, we have $G^+(z)=\log|\pi_1(z)|$; and for any $z\in V^-$, we have $G^-(z)=\log|\pi_2(z)|$. The sequence appearing inside of the limit on the right hand side is constant with respect to $n$, i.e., the limit is uniform.
\end{lmm}

\begin{proof}
For $(x,y)\in\A^{2,\an}_k(k)\cap V_R^+$, by the ultrametric inequality, we have
\begin{align*}
    |p(x)|=|x|^d>|x|\geq|y|.
\end{align*}
Since $d\geq 2$ and $R>|a|$, we also have $|p(x)|>|ay|$. Applying the ultrametric inequality again, we obtain
\begin{align*}
    |p(x)-ay|=|x|^d.
\end{align*}
Hence, we have $H(x,y)\in V^+_R$ again and $\|H(x,y)\|=|x|^d$, which gives
\begin{align*}
    \frac{1}{d}\log\|H(x,y)\|=\log|x|.
\end{align*}
As $H(x,y)\in V^+_R$, the same argument also gives
\begin{align*}
    \frac{1}{d^n}\log\|H^{\circ n}(x,y)\|=\log|x|
\end{align*}
for all $n=1,2,3,\ldots$. Then, since the set of $k$-rational points is dense in $\A^{2,\an}_k$, we have
\begin{align*}
    \frac{1}{d^n}\log\|H^{\circ n}(z)\|=\log|\pi_1(z)|
\end{align*}
for all $n$. Clearly the right hand side is independent of $n$, and the sequence $\{\frac{1}{d^n}\log\|H^{\circ n}(z)\|\}$ is thus constant with respect to $n$. The claim for $G^-(z)$ can be shown in the same way.
\end{proof}

To prove the claim for points outside $V^{\pm}_R$, we need the following properties, which also appears in the complex case.
\begin{prop} \label{NArelofV}
For any $R>\max\{|a_1|,\ldots,|a_d|,|a|,1\}$,
\begin{enumerate}
    \item $H(V^+_R)\subset V^+_R$, and $H(V^+_R\cup W_R)\subset V^+_R\cup W_R$,
    \item $H^{\circ -1}(V^-_R)\subset V^-_R$, and $H(V^-_R\cup W_R)\subset V^-_R\cup W_R$,
    \item $\A^{2,\an}\setminus K^+=\bigcup_{n\geq 0}H^{\circ-n}(V^+_R)$,
    \item $\A^{2,\an}\setminus K^-=\bigcup_{n\geq 0}H^{\circ n}(V^-_R)$.
    \end{enumerate}
\end{prop}

\begin{proof}

Note that the sets $K^+$ and $K^-$ are closed by the continuity of $H$. Therefore, it is sufficient to show the statements for classical points.

The first statement in (1) was already shown in the proof of Lemma \ref{lmm:NAtrivcase}. We show that $H(V^+_R\cup W_R)\subset V^+_R\cup W_R$,for which it is sufficient to show that $H(W_R)\subset V^+_R\cup W_R$. Taking any $(x,y)\in W_R$, since
$|\pi_2H(x,y)|=|x|\leq R$, we have $H(x,y)\in V^+_R$ if $|\pi_1H(x,y)|\geq R(\geq|\pi_2H(x,y)|)$ and $H(x,y)\in W_R$ if $|\pi_1H(x,y)|\leq R$. The same argument yields the second statement in (2), too.

For (3), note that $\|H(x,y)\|=|x|^d$ for any $(x,y)\in V^+_R(K)$. Hence, the forward orbit of any element in $V^+_R$ is unbounded, which shows that $\A^{2,\an}\setminus K^+\supset\bigcup_{n\geq 0}H^{\circ -n}(V^+_R)$. On the other had, we divide $\A^{2,\an}$ into three sets --- $V^+_R$, $V^-_R$ and $W_R$ --- and prove the statement for each of them. First, the case for $z\in\A^{2,\an}\cap V^+_R$ is trivial. Next, take $z\in\A^{2,\an}\cap W_R$. By (1), we have $H^{\circ n}(z)\in W_R$ for all $n$, or there exists some $N$ so that $H^{\circ n}(z)\in V^+_R$ for any $n\geq N$. As the orbit must be bounded in the former case, the latter case must hold, i.e., $x\in H^{\circ -N}(V^+_R)$. It remains to prove the statement for $z\in V^-_R\cap \A^{2,\an}_k\setminus \bigcup H^{\circ-n}(V^+_R)$. Assume the converse, i.e., $H^{\circ n}(z)\in V^-_R$, holds for any $n$. Then, since $|\pi_2H^{\circ n}(z)|=|\pi_1H^{\circ n-1}(z)|\leq|\pi_2H^{\circ n-1}(z)|$, the sequence $\{|\pi_2H^{\circ n}(z)|\}$ is non-increasing, i.e., bounded. As the orbit is unbounded, the set $\{|\pi_1 H^{\circ n}(z)|\}$ must be unbounded, which contradicts to the hypothesis on $V^-_R$, $|\pi_1 H^{\circ n}(z)|\leq|\pi_2H^{\circ n}(z)|$ for any $n$. A similar argument proves (4) by replacing $H$ with $H^{\circ-1}$.
\end{proof}

\begin{prop}
 The functions $G^+$ and $G^-$ are well-defined and continuous on $\A^{2,\an}$.
\end{prop}
\begin{proof}
The well-definedness easily follows from Proposition \ref{NArelofV}. Moreover, the continuity of $G^+$ in $\A^{2,\an}_k\setminus K^+$ was already shown in Lemma \ref{lmm:NAtrivcase}, and that in the interior of $K^+$ is easily seen since $G^+=0$ on $K^+$ by definition.

It thus remains only to show the continuity on $J^+$. Take any $z_0\in J^+$ and any $\varepsilon>0$. Suppose there exists a net $(z_{\lambda})_{\lambda\in\Lambda}$ in $\A^{2.an}_k$ such that $z_{\lambda}\to z_0$ and $G^+(z_{\lambda})\geq\varepsilon$ for all $\lambda$. Take $R$ and a natural number $N$ sufficiently large that $\sup|z_{\lambda}|<R$ and $R^d/d^N<\varepsilon$. Then, by setting $n(z)$ to be the smallest number such that $H^{\circ n}(z)\in V^+$, we have
\begin{align*}
    G^+(z)=\frac{1}{d^{n(z)}}\log\|H^{\circ n(z)}(z)\|.
\end{align*}
Also, by the ultrametric inequality, the absolute value of $H^{\circ n(z)}$ must be bounded above by $R^d$. Hence, for any $\lambda$, hence, we have $G^+(z_{\lambda})\leq R^d/d^{n(z)}$. Since we assume that $G^+(z_{\lambda})>\varepsilon>R^d/d^N$, we have $n(z_{\lambda})<N$, i.e., $H^{\circ n}(z_{\lambda})\in V^+_R$ for any $\lambda$. This is a contradiction by the assumption that $H^{\circ n}(z_0)\in W_R$ and the continuity of $H^{\circ n}$.

Therefore, for any $\varepsilon$, there exists an open neighborhood $U$ of $z_0$ such that $G^+(z)<\varepsilon$ for any $z\in U$; that is, the map $G^+$ is continuous around $z_0$. The continuity of $G^-$ can be shown in the same way.
\end{proof}

\begin{prop}\label{prop:thenonarchimedeangreenfunction}
 The sequences $\{G_n^+(z):=\frac{1}{d^n}\log^+\|H^{\circ n}(z)\|\}$ and $\{G_n^-(z):=\frac{1}{d^n}\log\|H^{\circ -n}(z)\|\}$ converge uniformly to $G^+$ and $G^-$, respectively, on each compact set of $\A^{2,\an}$.
\end{prop}

\begin{proof}
As noted at the beginning of this section, the family $\{W_R\}_R$ covers the whole $\A^{2,\an}_k$. We thus show the uniform convergence of these functions on $W_R$ for any $R$. For any $\varepsilon>0$, take $N$ sufficiently large that $(\log R)/d^{N-1}<\varepsilon$. Let $n(z)$ be the smallest $n$ such that $H^{\circ n}(z)\in W_R$. For $n\geq n(z)$, we have already shown that $G^+(z)=G^+_n(z)$, i.e., $|G^+(z)-G^+_n(z)|=0$. For $n(z)>n\geq N$, as $H^{\circ n(z)-1}(z)\in W_R$, we have 
\begin{align*}
0\leq G^+(z)=G_{n(z)}^+(z)=\frac{1}{d^{n(z)}}\log^+\|H(H^{\circ n(z)-1}(z)\|\leq\frac{1}{d^{n(z)}}\log R^d<\varepsilon.
\end{align*}
On the other hand, since $H^{\circ n}(z)\in W$,
\begin{align*}
    0\leq G^+_n(z)=\frac{1}{d^n}\log^+\|H^{\circ n}(z)\|<\frac{\log R}{d^{N-1}}<\varepsilon.
\end{align*}
Therefore, for any $n\geq N$, we have
\begin{align*}
    |G^+(z)-G_n^+(z)|<\varepsilon
\end{align*}
in any case.
\end{proof}

Clearly the functions $G^+$ and $G^-$ are continuous plurisubharmonic functions on $\A^{2,\an}_k$. Let their Laplacians be $T^{\pm}:=-\mathop{\mathrm{dd^c}}G^{\pm}$. From \cite{CLD}, the wedge product of two currents with continuous potentials is well-defined as in the complex pluri-potential theory; thus, we have an invariant measure $T^+\wedge T^-$.

Finally, since the function $G(z):=\max\{G^+(z),G^-(z)\}$ is also a continuous pluri-subharmonic function, we can define the following positive measure:
\begin{align*}
    \mu=\ddc{G}^{\wedge 2},
\end{align*}
which will play a central role later.

\begin{rmk}
In the complex dynamics of H\'enon maps, we have $\mu=T^+\wedge T^-$. It is natural to expect to have this equality to hold in the non-archimedean case, too, this remains an open question.
\end{rmk}

\section{A $1$-parameter family of complex H\'enon maps} \label{sec: Aoneparameterfamilyofcomplexhenonmaps}
In this section, we consider a (meromorphic) degenerating family of complex H\'enon maps over the unit punctured disk. The result presented here is basically similar to that in \cite{BS1}, but it is made somewhat more precise for application with the theory of hybrid spaces in the next section.

For $t\in\D^*$, we consider a complex H\'enon map
\begin{align}\label{eq:defoffamilyofhenonmaps}
    H_t(x,y)=(p_t(x)-a(t)y,x),
\end{align}
where the function $p_t\in\mathcal{O}(\mathbb{D})[t^{-1}][x]$ is a monic polynomial, and $a(t)\in\mathcal{O}(\mathbb{D})[t^{-1}]$ that does not vanish outside the origin. Recall that $\mathcal{O}(\mathbb{D})[t^{-1}]$ is the set of meromorphic functions on the open unit disk $\mathbb{D}$ with poles only at the origin. Note that, for any $t\neq0$, $H_t$ has the inverse $H_t^{\circ-1}(x,y)=(y, (p_t(y)-x)/a(t))$.

For each $t\neq0$, we have the Green functions $G_t^{\pm}(z)=\lim_{n}\frac{1}{d^n}\log^+\|H_t^{\circ n}(x)\|$ and the invariant measure $\mu_t=\ddc G_t^+\wedge \ddc G_t^-$. From \cite{BS1}, it is known that $\mu_t=(\ddc\max(G_t^+,G_t^-)^{\wedge2}$.

\begin{rmk}
Note that, to analyze $1$-parameter meromorphic degeneration of complex H\'enon maps, it is sufficient to consider the family of the above form. To see this, starting from the general form $H_t(x,y)=(a_0(t)x^d+\cdots+a_d(t)-a(t)y,b(t)x)$, where $a_0(t),\ldots,a_d(t),a(t),b(t)\in \mathcal{O}(\mathbb{D})[t^{-1}]$. Assume that $a_0(t)$ has a $(d-1)$-th root, i.e., an element $c(t)\in\mathcal{O}(\mathbb{D})[t^{-1}]$ such that $c(t)^d=a_0(t)$. This is enabled by considering the extension $\mathcal{O}(\mathbb{D})[t^{-1}]$ of a larger ring $\mathcal{O}(\mathbb{D})[t^{-1/(d-1)}]$, or, more precisely, writing
\begin{align*}
    a_0(t)&=c^{(m)}t^m+c^{(m+1)}t^{m+1}+\cdots \\
    &=c^{(m)}t^m(1+(c^{(m+1)}t+c^{(m+2)}t^2+\cdots)).
\end{align*}
By the binomial theorem, we have a $(d-1)$-th root of $1+(c^{(m+1)}t+c^{(m+2)}t^2+\cdots)$, denoted by $f(t)$. Then, setting $t=s^{d-1}$, we have a $(d-1)$-th root of $a_0$, i.e., $c(s)=(c^{(m)})^{1/d-1}s^m\times f(s^{d-1})$. Take an affine transformation $\sigma$ of $\C^2$ such that
\begin{align*}
    \sigma_s(x,y)=(c(s)x,b(s^{d-1})c(s)y).
\end{align*}
Then, we obtain a new H\'enon map that is conjugate to $H_{s^{d-1}}$:
\begin{align*}
    \sigma_s\circ H_{s^{d-1}}\circ\sigma_s^{-1}(x,y)=(x^d+a'_1(s)x^{d-1}+\cdots+a'_d(s)+a(s)b(s^{d-1})y.x),
\end{align*}
where $a'_i(s)$ is a multiple of $c(s)$, and $a_i(s^{d-1})$ for $i=1,\ldots,d$.
\end{rmk}

From this point, we discuss families of complex H\'enon maps of the form given in equation (\ref{eq:defoffamilyofhenonmaps}). As before, we define the functions $G_t$, with the addition of several notations fo later argument. First, define the constant $R_t$ as
\begin{align*}
    R_t:=c\max\{1,|a(t)|,|a(t)|^{-1},|a(t)|^{-2},|a_1(t)|,\ldots,|a_d(t)|\}.
\end{align*}
The constant $c>1$ is an absolute constant (i.e., it is chosen independently from $t$) that will be determined in later (it will turn out that $c=5$ suffices). Also, as in Section \ref{sec:non-archimedean-henon-maps}, certain special subsets of $\C^2$ are defined:
\begin{align*}
    &V_t^+=\left\{(x,y)\in\C^2\middle| |x|\geq|y|\text{, and }|x|\geq R_t\right\}, \\
    &V_t^-=\left\{(x,y)\in\C^2\middle| |y|\geq|x|\text{, and }|y|\geq R_t\right\}, \\
    &W_t=\left\{(x,y)\in\C^2\middle| |x|\leq R_t\text{, and }|y|\leq R_t\right\}, \\
    &K_t^+=\left\{z\in\C^2\ \middle|\ \{\|H_t^n(z)\|\}_n\text{ is bounded}\right\}, \\
    &K_t^+=\left\{z\in\C^2\ \middle|\ \{\|H_t^{-n}(z)\|\}_n\text{ is bounded}\right\}.
\end{align*}
Throughout this section, the norm on $\C^2$ is taken as $\|(x,y)\|:=\max\{|x|,|y|\}$. 
Let $G_{n,t}^{\pm}$ be as in Section \ref{sec:non-archimedean-henon-maps}; that is,
\begin{align*}
    G_{n,t}^{\pm}(x,y)&:=\frac{1}{d^n}\log^+\|H^{\circ n}(x,y)\|, \\
    G_t^{\pm}(x,y)&:=\lim_{n\to\infty}G_{n,t}^{\pm}(x,y).
\end{align*}
Lastly, define $G_t:=\max\{G_t^+,G_t^-\}$ and $G_{n,t}:=\max\{G_{n,t}^+,G_{n,t}^-\}$, and for later arguments, define
\begin{align*}
    (x_n,y_n):=H^{\circ n}(x,y).
\end{align*}

The goal of this section is the following estimation.
\begin{prop}\label{prop:mainevaluation}
 Let $\{H_t\}_{t\in\D^*}$ be a family of H\'enon maps as in equation (\ref{eq:defoffamilyofhenonmaps}). There exists a decreasing sequence of positive real numbers $\{\epsilon_n\}_{n=2}^{\infty}$ converging to $0$ such that
\begin{align*}
    |G_t(z)-G_{n,t}(z)|\leq\epsilon_n\log|t|^{-1}
\end{align*}
for all $n=2,3,4,\ldots$ and $t$ sufficiently close to $0$.
\end{prop}

Before the estimation, several lemmas are necessary.
\begin{lmm}\label{lemma:H(x,y)/(x,y)}
  There exists a positive real number $\delta$, chosen independently of a parameter $t$, such that for any $t$ and $(x,y)\in V^+_t$, we have
  \begin{align*}
    (1-\delta)|x|^d\leq\|H_t(x,y)\|\leq(1+\delta)|x|^d.
  \end{align*}
  Also, for any $t$ and $(x,y)\in V^-_t$,
  \begin{align*}
    |a(t)|^{-1}(1-\delta)|y|^d\leq\|H_t^{\circ -1}(x,y)\|\leq|a(t)|^{-1}(1+\delta)|y|^d.
  \end{align*}

The constant $\delta$ can be written explicitly:
\begin{align*}
    \delta=\frac{c^d+c^2-c-1}{c^{d+1}-c^d}.
\end{align*}
\end{lmm}

 Note that even though the constant $\delta$ is independent of $t$, the definition of the sets $V_t^{\pm}$ depends on $t$. This lemma was originally proven in \cite{BS1}, but we give a detailed proof here to explicitly indicate the constant $\delta$ for the given $R_t$, which is essential in the proof of Proposition \ref{prop:mainevaluation}.
  \begin{proof}
  For any $(x,y)\in V_t^+$, we have
  \begin{align*}
      \|H_t(x,y)\|&=|x_1| \\
      &=|x^d+a_1(t)x^{d-1}+\cdots+a_d(t)+a(t)y| \\
      &=|x|^d\left|1+\frac{a_1(t)}{x}+\cdots+\frac{a_d(t)}{x^d}-\frac{a(t)y}{x^d}\right| \\
      &\leq|x|^d\left(1+\left|\frac{a_1(t)}{x}\right|+\cdots+\left|\frac{a_d(t)}{x^d}\right|-\left|\frac{a(t)y}{x^d}\right|\right).
  \end{align*}

Similarly,
  \begin{align*}
      \|H_t(x,y)\|&=|x_1|\geq|x|^d\left(1-\left(\left|\frac{a_1(t)}{x}\right|+\cdots\left|\frac{a_d(t)}{x^d}\right|+\left|\frac{a(t)y}{x^d}\right|\right)\right).      
  \end{align*}
  We seek an upper bound on the terms $|a_1(t)/x|+\cdots+|a_d(t)/x^d|+|ay/x|$. By the definitions of $R_t$ and $V_t^+$, we have
  \begin{align*}
      &|x|\geq R\geq c, \\
      &|a_i(t)/x|\leq c\text{ and } |a(t)/x|\leq c\text{ for all }i,\text{ and} \\
      &|y|\leq|x|.
  \end{align*}
  Hence, $|a_i(t)/x^i|\leq c^{-i}$ and $|a(t)y/x^d|\leq c^{-(d-1)}$. Then, we have
  \begin{align*}
      \left|\frac{a_1(t)}{x}\right|+\cdots\left|\frac{a_d(t)}{x^d}\right|+\left|\frac{a(t)y}{x^d}\right|\leq \frac{1}{c}+\cdots\frac{1}{c^d}+\frac{1}{c^{d-1}} =\frac{c^d+c-1}{c^{d+1}-c^d}.
  \end{align*}

Similarly, for $(x,y)\in V_t^-$,
\begin{align*}
      \|H_t^{\circ -1}(x,y)\|&=|y_{-1}| \leq\frac{|y|^d}{|a(t)|}\left(1+\left|\frac{a_1(t)}{y}\right|+\cdots\left|\frac{a_d(t)}{y^d}\right|+\left|\frac{x}{y^d}\right|\right),
\end{align*}
and
\begin{align*}
      \|H_t^{\circ -1}(x,y)\|=|y_{-1}| \geq\frac{|y|^d}{|a(t)|}\left(1-\left(\left|\frac{a_1(t)}{y}\right|+\cdots\left|\frac{a_d(t)}{y^d}\right|+\left|\frac{x}{y^d}\right|\right)\right),
\end{align*}
and
  \begin{align*}
      \left|\frac{a_1(t)}{y}\right|+\cdots\left|\frac{a_d(t)}{y^d}\right|+\left|\frac{x}{y^d}\right|&\leq \frac{c^d+c^2-c-1}{c^{d+1}-c^d}.
  \end{align*}
Thus, taking $\delta=(c^d+c^2-c-1)/(c^{d+1}-c^d)$ proves the claim.
  \end{proof}

\begin{lmm}\label{lmm: G_t^+forV_t^+}
For any $(x,y)\in V^+_t$ and $n\in\N$,
\begin{align*}
    |G_t^+(x,y)-G_{n,t}^+(x,y)|\leq\left(\sum_{i>n}\frac{1}{d^i}\right)\log(1-\delta)^{-1}.
\end{align*}
Similarly, for any $z\in V^-_t$ and $n\in\N$,
\begin{align*}
    |G_t^-(x,y)-G_{n,t}^-(x,y)|\leq\left(\sum_{i>n}\frac{1}{d^i}\right)\log|a(t)|^{-1}(1-\delta)^{-1}.
\end{align*}
\end{lmm}

\begin{proof}
By Lemma \ref{lemma:H(x,y)/(x,y)}, for any $(x,y)\in V_t^+$,
\begin{align*}
    |G_{n+1,t}^+(x,y)-G_{n,t}^+(x,y)|&=\left|\frac{1}{d^{n+1}}\log^+\|H_t^{\circ n+1}(x,y)\|-\frac{1}{d^{n}}\log^+\|H_t^{\circ n}(x,y)\|\right| \\
    &=\frac{1}{d^{n+1}}\left|\log^+\frac{\|H_t(H_t^{\circ n}(x,y))\|}{\|H_t^{\circ n}(x,y)\|^d}\right| \\
    &\leq \frac{1}{d^{n+1}}\max\{|\log(1-\delta)|,|\log(1+\delta)|\} \\
    &=\frac{1}{d^{n+1}}\log(1-\delta)^{-1},
\end{align*}
where the last equality follows from $(1-\delta)^{-1}\geq1+\delta$, given $\delta<1$. The lemma then follows since $|G_t^+-G_{n,t}^+|\leq\sum_{i\geq n}|G_{n+1,t}^+-G_{n,t}^+|$. A similar argument proves the other assertion.
\end{proof}
  
 Next, we need to consider $G_t^+$ (resp.\ $G_t^-$) on $V_t^-$ (resp.\ $V_t^+$). In this case, in general, it is difficult to estimate the speed of convergence. However, since the goal is to estimate $G_t$ which is the maximum of $G_t^+$ and $G_t^-$, the speed of convergence of each $G_t^+$ and $G_t^-$ is of no concern according to the following claim.
  
  \begin{lmm}\label{lmm: n^+(z)}
  For $z\in\C^2$, let $n^+(z)$ (resp.\ $n^-(z)$) be the smallest number $n$ such that $H_t^{\circ n}(z)\in V_t^+$ (resp.\ $H_t^{\circ -n}(x,y)\in V_t^-$). If $H_t^{\circ n}(z)\not\in V_t^+$ (resp. $H_t^{\circ -n}(z)\not\in V_t^-$) for any $n\in\N$ , we set $n^+(z)=\infty$ (resp. $n^-(z)=\infty$). Then, for sufficiently large $c$ and any $z\in V_t^+$, either of the following conditions holds:
  \begin{enumerate}
      \item $n^-(z)=1$;
      \item $G_t^+(z)\geq G_t^-(z)$.
  \end{enumerate}
   Similarly, for any $z\in V_t^-$, either of the following conditions holds:
  \begin{enumerate}
      \item $n^+(z)=1$;
      \item $G_t^-(z)\geq G_t^+(z)$.
  \end{enumerate}
  \end{lmm}
  
\begin{proof}
First, We claim that for any $n<n^-(z)$ and $z=(x,y)\in V_t^+$, 
\begin{align}
    \|H_t^{\circ -n}(x,y)\|\leq\|(x,y)\|(=|x|). \label{eq:lemma3.5-1}
\end{align}

Indeed, if $H_t^{\circ -1}(x,y)\in V_t^+$, then $|x_{-1}|=|y|\leq|x|$ and $|y_{-1}|\leq|x_{-1}|\leq|x|$ by definition. If $H_t^{\circ -1}(x,y)\in W_t$, then $|x_{-1}|,$ $|y|\leq R_t\leq|x|$. If $n^-(z)\geq n>1$, there are two cases, $H_t^{\circ -1}(x,y)\in V_t^+$ and $H_t^{\circ -1}(x,y)\in W_t$. In the former case of $H_t^{\circ -1}(x,y)\in V_t^+$, the same argument yields $\|H_t^{\circ -2}(z)\|\leq|x_{-1}|\leq|x|$. In the latter case of $H_t^{\circ -1}(x,y)\in W_t$, since $H_t^{\circ -1}(W_t)\subset W_t\cup V_t^-$, we must have $H_t^{\circ -2}(z)\in W_t$. Hence, $|x_{-2}|$, $|y_{-2}|\leq R_t\leq|x|$. We can apply the same argument for $n<n^-(z)$.

Next, we prove the following two inequalities:
\begin{align}
    G_t^-(z)&\leq\frac{1}{d^{n^-(z)-1}}\log\|z\|+\frac{1}{d^{n^-(z)-1}}\log(1-\delta)^{-1}+\frac{1}{d^{n^-(z)}}\frac{1}{d-1}\log|a(t)|^{-1}(1-\delta)^{-1}, \label{eq:estimateG_t^-} \\
    G_t^+(z)&\geq\log\|z\|-\frac{1}{d-1}\log(1-\delta)^{-1}. \label{eq:estimateG_t^+}
\end{align}

First, rearding the inequality (\ref{eq:estimateG_t^-}), for $z=(x,y)\in V^+_t$, since $H_t^{-n^-(z)}\in V_t^-$, we have
\begin{align*}
    |G_t^-(z)-G_{n^-(z),t}^-(z)|&=\frac{1}{d^{n^-(z)}}|G_t^-(H_t^{-n^-(z)}(z))-(\log^+\|H_t^{-n^-(z)}(z)\|)| \\
    &\leq\frac{1}{d^{n^-(z)}}\left(\sum_{i>0}\frac{1}{d^i}\right)\log|a(t)|^{-1}(1-\delta)^{-1}\\
    &=\frac{1}{d^{n^-(z)}}\frac{1}{d-1}\log|a(t)|^{-1}(1-\delta)^{-1}.
\end{align*}
By inequality (\ref{eq:lemma3.5-1}), for any $n<n^-(z)$, we have 
\begin{align*}
    \log\|H_t^{-n}(z)\|\leq\log\|z\|.
\end{align*}
Next, for $n=n^-(z)$, by  Lemma \ref{lemma:H(x,y)/(x,y)},
\begin{align*}
    G_{n^-(z),t}^-(z)&=\frac{1}{d^{n^-(z)}}\log\|H_t^{\circ-n^-(z)}(z)\| \\
    &=\frac{1}{d^{n^-(z)}}\log\|H_t^{\circ-1}(H_t^{\circ -(n^-(z)-1)}(z)\| \\
    &=\frac{1}{d^{n^-(z)-1}}\log\|z\|+\frac{1}{d^{n^-(z)}}\log|a(t)|^{-1}(1-\delta)^{-1}.
\end{align*}
By combining the above results, we obtain
\begin{align*}
    G_t^-(z)&\leq G_{n^-(z),t}^-(z)+\frac{1}{d^{n^-(z)}}\frac{1}{d-1}\log|a|^{-1}(1-\delta)^{-1} \\
    &\leq\frac{1}{d^{n^-(z)-1}}\log\|z\|+\frac{1}{d^{n^-(z)-1}}\log(1-\delta)^{-1}+\frac{1}{d^{n^-(z)}}\frac{1}{d-1}\log|a|^{-1}(1-\delta)^{-1}.
\end{align*}

Next, for inequality (\ref{eq:estimateG_t^+}), Since $z\in V_t^+$, we have
\begin{align*}
    |G_t^+(z)-\log\|z\||&=|G_t^+(z)-G_{0,t}^+(z)|\\
    &\leq\left(\sum_{i>0}\frac{1}{d^i}\right)\log(1-\delta)^{-1}\\
    &=\frac{1}{d-1}\log(1-\delta)^{-1}.
\end{align*}
Hence,
\begin{align*}
    G_t^+(z)\geq\log\|z\|-\frac{1}{d-1}\log(1-\delta)^{-1}.
\end{align*}
It is then sufficient to show that
\begin{align*}
    \log\|z\|-\frac{1}{d-1}\log(1-\delta)^{-1}\geq\frac{1}{d^{n^-(z)-1}}\log\|z\|+\frac{1}{d^{n^-(z)-1}}\log(1-\delta)^{-1}+\frac{1}{d^{n^-(z)}}\frac{1}{d-1}\log|a(t)|^{-1}(1-\delta)^{-1},
\end{align*}
which is equivalent to
\begin{align*}
    \left(1-\frac{1}{d^{n^-(z)-1}}\log\|z\|\right)\geq\frac{1}{d^{n^-(z)}}\left(1+\frac{1}{d-1}\right)\log|a(t)|^{-1}+\left(\frac{1}{d^{n^-(z)}}\left(1+\frac{1}{d-1}\right)+\frac{1}{d-1}\right)\log(1-\delta)^{-1}.
\end{align*}

Note that this inequality holds for sufficiently large $c$ in the definition of $R_t$ if $|a(t)|$ is bounded by below. Consider the case of $a(t)\to0$ as $t\to0$. It is sufficient to prove the inequality when the left hand side takes its maximum and the right hand side takes its minimum, which occurs when $d=2$ and $n^-(z)=2$ on both hand sides. The inequality thus becomes
\begin{align*}
    \frac{1}{2}\log\|z\|\geq \frac{1}{2}\log|a(t)|^{-1}+\frac{3}{2}\log(1-\delta)^{-1},
\end{align*}
or equivalently,
\begin{align*}
    \|z\|\geq\frac{1}{(1-\delta)^{-3}}|a(t)|^{-1}.
\end{align*}
As $\|z\|\geq c\max\{|a(t)|,|a(t)|^{-1},|a(t)|^{-1},|a_1(t)|,\ldots,|a_d(t)|\}\geq c|a(t)|^{-1}$, this inequality holds when
\begin{align*}
    c>(1-\delta)^{-3}=\left(1-\frac{c^d+c^2-c-1}{c^{d+1}-c^d}\right)^{-3}.
\end{align*}
Since the term $1-(c^d+c^2-c-1)/(c^{d+1}-c^d)$ tends to $1$ as $c\to\infty$, this inequality is always true for sufficiently large $c$, and numerical calculation demonstrates that $c=5$ is sufficient.

A similar argument proves the inequality for $z\in V_t^-$, and the function $G_t^+$ can be estimated from above as
\begin{align*}
    G_t^+(z)\leq \frac{1}{d^{n^+(z)-1}}\log\|z\|+\left(\frac{1}{d^{n^+(z)}}+\frac{1}{d^{n^+(z)}(d-1)}\right)\log(1-\delta)^{-1}.
\end{align*}
Next, $G_t^-(z)$ can be estimated from below as
\begin{align*}
    G_t^-(z)\geq\log\|z\|-\frac{1}{d-1}\log|a(t)|^{-1}(1-\delta)^{-1}.
\end{align*}
By combining the above results, it is sufficient to show that
\begin{align*}
    \frac{1}{d^{n^+(z)-1}}\log\|z\|+\left(\frac{1}{d^{n^+(z)}}+\frac{1}{d^{n^+(z)}(d-1)}\right)\log(1-\delta)^{-1}\leq \log\|z\|-\frac{1}{d-1}\log|a(t)|^{-1}(1-\delta)^{-1}.
\end{align*}
Setting $n^+(z)=2$ and $d=2$, we have
\begin{align*}
    \frac{1}{2}\log\|z\|\geq \log|a(t)|^{-1}+\frac{3}{2}\log(1-\delta)^{-1}.
\end{align*}
Finally since $\|z\|\geq c|a(t)|^{-2}$ by definition, this inequality holds for the same $c$ chosen above.
\end{proof}

The rest of this section constructs the proof of Proposition \ref{prop:mainevaluation} by dividing $\C^2$ into $V_t^+$, $V_t^-$, and $W_t$. The notations $n^+(z)$ and $n^-(z)$ are used as in the proof of Lemma \ref{lmm: n^+(z)}.

\subsection{Case 1: when $z\in W_t$}

Since $|G_t-G_{n,t}|\leq\max\{|G_t^+-G_{n,t}^+|,|G_t^--G_{n,t}^-|\}$, it is sufficient to show that $|G_t^+-G_{n,t}^+|\leq\epsilon_n\log|t|^{-1}$ and $|G_t-G_{n,t}^-|\leq\epsilon_n\log|t|^{-1}$. The proofs of these inequalities are similar, and we only present the estimate for $|G_t^+-G_{n,t}^+|$.

As $G_{n,t}^+\geq0$ and $G_t^+\geq0$ by definition, we have
\begin{align*}
    |G_t^+-G_{n,t}^+|\leq G_t^++G_{n,t}^+.
\end{align*}
Hence, it is sufficient to estimate $G_t^+$ and $G_{n,t}^+$ from above.

For any $z\in W_t$ and any $n\geq n^+(z)$, we can apply Lemma \ref{lmm: G_t^+forV_t^+} as follows:
\begin{align*}
    |G_t^+(z)-G_{n,t}^+(z)|&=\frac{1}{d^{n^+(z)}}\left|G_t^+(H^{\circ n^+(z)}(z))-\frac{1}{d^{n-n^+(z)}}\log\|H_t^{\circ n-n^+(z)}(H_t^{n^+(z)}(z))\|\right| \\
    &=\frac{1}{d^{n^+(z)}}\left|G_t^+(H^{n^+(z)}(z))-g_{n-n^+(z),t}(H^{n^+(z)}(z))\right| \\
    &\leq \frac{1}{d^n(d-1)}\log(1-\delta)^{-1}.
\end{align*}
 
 This tends to be $0$ as $n\to\infty$, and the speed of convergence is independent of $t$.
 
 For $n<n^+(z)$, we have $H_t^{\circ n}(z)\in W_t$ since $H_t(W_t)\subset W_t\cup V_t^+$, i.e., $\|H_t^{\circ n}(z)\|\leq R_t$. Hence,
 \begin{align*}
     G_{n,t}^+(z)\leq \frac{1}{d^n}\log R_t.
 \end{align*}
 On the other hand, we can also estimate $G_t^+$ as in the proof of Lemma \ref{lmm: n^+(z)}. As $\|H_t^{n^+(z)-1}(z)\|\leq R_t$, we have 
 \begin{align*}
     \log\|H_t^{\circ n^+(z)}(z)\|\leq d\log R_t+\log(1-\delta)^{-1}.
 \end{align*}
 Then, since $H_t^{n^+(z)}(z)\in V_t^+$, application of Lemma \ref{lmm: G_t^+forV_t^+} yields
 \begin{align*}
     G_t^+(z)&\leq \frac{1}{d^{n^+(z)-1}}\log R_t+\frac{d}{d^{n^+(z)}(d-1)}\log(1-\delta)^{-1} \\
     &\leq \frac{1}{d^{n-1}}\log R_t+\frac{1}{d^{n-1}(d-1)}\log(1-\delta)^{-1}.
 \end{align*}
 Hence, for $n<n^+(z)$, we have 
 \begin{align*}
     |G_t^+-G_{n,t}^+|\leq\frac{1}{d^{n-1}}\log R_t+\frac{1}{d^{n-1}(d-1)}\log(1-\delta)^{-1}.
 \end{align*}
 By the definition of $R_t$, there exists a real number $\alpha>0$ such that
 \begin{align*}
     \frac{1}{d^{n-1}}\log R_t+\frac{1}{d^{n-1}(d-1)}\log(1-\delta)^{-1}\leq \frac{\alpha}{d^n}\log|t|^{-1}.
 \end{align*}
Therefore, for all $n$ larger or smaller than $n^+(z)$, we have
\begin{align*}
    |G_t^+(z)-G_{n,t}^+(z)|\leq\frac{1}{d^{n-1}}\log R_t+\frac{1}{d^{n-1}(d-1)}\log(1-\delta)^{-1}\leq \frac{\alpha}{d^n}\log|t|^{-1}.
\end{align*}

The same argument can be applied for $|G_t^--G_{n,t}^-|$. For $n\geq n^-(z)$,
\begin{align*}
    |G_t^-(z)-G_{n,t}^-(z)|\leq\frac{1}{d^n(d-1)}\log|a(t)|^{-1}(1-\delta)^,
\end{align*}
and for $n<n^-(z)$,
\begin{align*}
    |G_t^-(z)-G_{n,t}^-(z)|&\leq\frac{1}{d^{n-1}}\log R_t+\frac{1}{d^{n-1}(d-1)}\log|a(t)|^{-1}(1-\delta)^{-1} \\
    &\leq \frac{\alpha}{d^n}\log|t|^{-1}+\frac{1}{d^{n-1}(d-1)}\log|a(t)|^{-1}
\end{align*}
By definition, there exists a real number $\beta>0$ such that
\begin{align*}
    \frac{1}{d^{n-1}(d-1)}\log|a(t)|^{-1}\leq\frac{\beta}{d^n}\log|t|^{-1}.
\end{align*}
Then, we have
\begin{align*}
    |G_t^-(z)-G_{n,t}^-(z)|\leq\frac{\alpha+\beta}{d^n}\log|t|^{-1},
\end{align*}
and therefore,
\begin{align}\label{eq:G_t-G_n,tforW_t}
    |G_t-G_{n,t}|\leq\frac{\alpha+\beta}{d^n}\log|t|^{-1}.
\end{align}

\subsection{Case 2: when $z\in V_t^+$}

Since the same estimation as in Case 1 is valid for $n\geq n^-(z)$ (note that $n^+(z)=0$), we have
\begin{align}\label{eq:G_t-G_n,tforV_t^+}
    |G_t-G_{n,t}|\leq \frac{\alpha+\beta}{d^n}\log|t|^{-1}.
\end{align}
for all $n\geq n^-(z)$.

For $2\leq n<n^-(z)$, we prove the two inequalities $G^+_t(z)\geq G^-_t(z)$ and $G_{n,t}^+(z)\geq G_{n,t}^-(z)$, which together imply that
\begin{align*}
    |G_t(z)-G_{n,t}(z)|=|G_t^+(z)-G_{n,t}^+(z)|\leq\frac{\alpha+\beta}{d^n}\log|t|^{-1}.
\end{align*}
While the inequality $G^+_t(z)\geq G^-_t(z)$ follows from Lemma \ref{lmm: n^+(z)}, it remains to show that $G_{n,t}^+(z)\geq G_{n,t}^-(z)$. First,
\begin{align*}
    G_{n,t}^+(z)\geq\log\|z\|-\frac{1-\frac{1}{d^n}}{d-1}\log(1-\delta)^{-1}.
\end{align*}
As $z\in V_t^+$ and $n<n^-(z)$, we have $\|z\|\geq\|H_t^{\circ -n}(z)\|$. Hence,
\begin{align*}
    G_t^-(z)\leq\frac{1}{d^n}\log\|z\|.
\end{align*}
The inequality $G_{n,t}^+\geq G_{n,t}^-$ follows from these estimations; indeed, it is sufficient to show that
\begin{align*}
    \log\|z\|-\frac{1-\frac{1}{d^n}}{d-1}\log(1-\delta)^{-1}\geq \frac{1}{d^n}\log\|z\|,
\end{align*}
which is equivalent to
\begin{align*}
    \log\|z\|\geq \frac{1}{d-1}\log(1-\delta)^{-1}.
\end{align*}
Finally, following the proof of Lemma \ref{lmm: n^+(z)}, we determine $c$ such that
\begin{align*}
    c>(1-\delta)^{-3}.
\end{align*}
For such $c$, we can easily show that $c>(1-\delta)^{-1/(d-1)}$ also holds. Since $\|z\|\geq R_t\geq C$, this proves the inequality.

\subsection{Case 3: when $z\in V_t^-$}

The basic argument is the same as in Case 2. For $n\geq n^+(z)$,
\begin{align}\label{eq:G_t-G_n,tforV_t^-}
        |G_t-G_{n,t}|\leq\frac{\alpha+\beta}{d^n}\log|t|^{-1}.
\end{align}
Assuming $n^+(z)\geq 2$, by applying the same argument as in Case $2$, we have $G_t^+\leq G_t^-$ and $G_{n,t}^+\leq G_{n,t}^-$ for all $n<n^+(z)$. Hence, the inequality (\ref{eq:G_t-G_n,tforV_t^-}) actually holds for all $n\geq2$.

Therefore, by inequalities (\ref{eq:G_t-G_n,tforW_t}), (\ref{eq:G_t-G_n,tforV_t^+}) and (\ref{eq:G_t-G_n,tforV_t^-}), for any $z\in \C^2$, we have
\begin{align}\label{eq:G_t-G_n,t}
    |G_t-G_{n,t}|\leq\frac{\alpha+\beta}{d^n}\log|t|^{-1}.
\end{align}



\section{Monge-Amp\`ere measures on the hybrid projective plane} \label{sec:mongeamperemeasuresonthehybridprojectiveplane}

In this section, we adapt certain basic notions and facts on hybrid dynamics by Favre \cite{Fa} adapted to our particular situation, for application to the family of H\'enon maps in the next section.

For a trivial family $\pi:\pr^2(\C)\times\D^*\to\D^*$ of $\pr^2(\C)$ over a punctured disk $\D^*$, we can take its trivial snc model as $\tilde{\pi}:\pr^2(\C)\times\D\to\D$ of $\pi$. During this section we fix a coordinate $([X:Y:Z],t)\in\pr^2(\C)\times\D$, which is necessary to fix a Fubini-Study metric later.

Let $\mathcal{L}$ be the pull-back of $\Sh(1)$ by the projection $\pr^2(\C)\times\D\to\pr^2(\C)$. A datum $\F:=\{d,\tau_1,\ldots,\tau_n,\mathfrak{a},p_{\mathcal{X}}\}$ is a regular admissible datum if $\F$ consists of the following:
\begin{itemize}
    \item A non-zero natural number $d\in\N\setminus\{0\}$,called the degree of $\F$;
    \item A  vertical fractional ideal $\mathfrak{a}$ satisfying $t^N\mathfrak{a}\subset\mathcal{O}_{\pr^2(\C)\times\D}$;
    \item Meromorphic sections $\tau_1,\ldots,\tau_n$ of $\mathcal{L}^{\otimes d}$ generating $\mathfrak{a}$ which is holomorphic off to the central fibre;
    \item A log-resolution $p_{\mathcal{X}}:\mathcal{X}\to\pr^2(\C)\times\D$ of $\mathfrak{a}$.
\end{itemize}

Here, the log-resolution of $\mathfrak{a}$ is an snc model $p:\mathcal{X}\to\pr^{2}(\C)\times\D$ such that $p^*\mathfrak{a}$ is equal to $\mathcal{O}_{\mathcal{X}(D)}$ for some vertical divisor $D$.
For a regular admissible datum $\F$, we next define the model functions and highlight several facts for use in later arguments.


\subsection{Complex Monge-Amp\`ere measure} \label{subsec:complexmongeamperemeasure}
Over the complex variety $\pr^2(\C)\times\D^*$, a model function $\phi_{\F}$ is attached to each regular admissible datum $\F$ by
\begin{align*}
    \phi_{\F}(z):=\log\max\{|\tau_1(z)|_*,\ldots,|\tau_n(z)|_*\},
\end{align*}
where the metric $|\cdot|_*$ is 
the hermitian metric given by $(|X|^2+|Y|^2+|Z|^2)^{-1}$, with $[X:Y:Z]$ the homogeneous coordinate of the projective plane that was fixed in the beginning of this section.

To apply the theory to a dynamical situation, the uniformity of functions plays a central role.
\begin{dfn}
A function $\phi:\pr^2(\C)\times\D^*\to\R$ is \textit{uniform} if there exists a sequence of model functions $\{\phi_{\F_n}\}$ such that
\begin{itemize}
    \item the degree $d_n$ of $\F_n$ diverges as $n\to\infty$, and
    \item there exists a decreasing sequence $\{\epsilon_n\}_n$ of positive real numbers converging to $0$ (which is independently of $t$) such that for any $t\in\D^*$,
    \begin{align*}
        \sup_{z\in \pr^2(\C)}\left|\phi(z,t)-\frac{1}{d_n}\phi_{\F_n}(z,t)\right|\leq \epsilon_n \log|t|^{-1}.
    \end{align*}
\end{itemize}
\end{dfn}

Note that a uniform function defines a continuous metric on $\mathcal{L}|_{\pr^2(\C)\times\D^*}$ as follows:
\begin{align}
    |\cdot|_{\phi}:=|\cdot|_*e^{-\phi}. \label{eq:complexmetricinducedbyfunctions}
\end{align}

On each $t\in\D^*$, let $|\cdot|_{\phi(\cdot,t)}$ be the metric restricting $|\cdot|_{\phi}$ to the fiber $\pr^2(\C)\times\{t\}$. The Monge-Amp\`ere measure $\MA_t(\phi)$ of $\phi(\cdot,t)$ is then defined as $(\ddc\log|\cdot|_{\phi(\cdot,t)})^{\wedge 2}/4\pi^2$.


\subsection{Non-archimedean Monge-Amp\`ere measure} \label{section:non-archimedeanmongeamperemeasure}
This definition is highly similar to that given above for the complex Monge-Amp\`ere measure. For a regular admissible datum $\F$, define the model function $g_{\F}:\pr^{2,an}_{\C((t))}\to\R$ as follows. First, take an open coverings $U_1:=\{X\neq0\}$, $U_2:=\{Y\neq0\}$, and $U_3:=\{Z\neq0\}$ of $\pr^2(\C)\times\D$. By Lemma 2.6 in \cite{Fa} we can take a set of meromorphic functions $\{f^i, g^i_1,\ldots,g^i_l\}$ on $U_i$ such that $f^i\mathfrak{a}|_{U_i}$ is genenerated by $\{g^i_1,\ldots,g^i_l\}$. Then, we may regard $f^i$ and $g^i_j$'s as elements of a polynomial ring with coefficients in $\C((t))$. Taking a base change of the $U_i$ to $\Spec\C((t))$ and the Berkovich analytification, we obtain open coverings $U_{i,\C((t))}^{\an}$ of $\pr^{2,\an}_{\C((t))}$. For any $x\in U_{i,\C((t))}^{\an}$, define 
\begin{align*}
    g_{\F}(x):=\log|\mathfrak{a}|(x)=\inf_j\{\log|g^i_j(x)|\}-\log|f^i(x)|.
\end{align*}
This definition is independent of the choices of the pen coverings and the generators $f^i$ and $g_j^i$. 

Next, define the canonical metric on $\mathcal{L}_{\C((t))}^{\an}$ over $\pr^{2,\an}_{\C((t))}$. For any $z\in\pr^{2,\an}_{\C((t))}$, we first assume $z\in\{Z\neq 0\}$ and take an inhomogeneous coordinate $(x(z),y(z))$, where $x(z)$ and $y(z)$ are regarded as elements in the complete residue field $\mathcal{H}(z)$. Then, we can write the point $z$ via a homogeneous coordinate in $\pr^{2,\an}_{\mathcal{H}(z)}$, e.g., $z=[x(z):y(z):1]$. For any homogeneous polynomial $P$, we can consider $P(z)=P([x(z):y(z):1])\in\mathcal{H}(z)$ and its norm $|P(z)|$, since the field $\mathcal{H}(z)$ has a natural norm induced from the seminorm $z$. If the point $z$ belongs to the locus $\{Z=0\}$, then we apply the same argument by replacing $Z$ with $X$ or $Y$. For any global section $\tau$ of $\mathcal{L}^{\an}_{\C((t))}$, define
\begin{align*}
    |\tau(z)|_{\can}:=\frac{|P_{\tau}(z)|}{\max\{|X(z)|,|Y(z)|,|Z(z)|\}},
\end{align*}
where $P_{\tau}$ is a homogeneous polynomial that is attached to $\tau$, and must be of degree $1$. Since it entails the quotient of the norms of two homogeneous polynomials of the same degree, the value $|\tau(z)|_{\can}$ is independent of the choice of homogeneous coordinates. The function $z\mapsto |\tau(z)|_{\can}$ is continuous by definition.

The metric $|\cdot|_{\can}$ naturally induces a metric $|\cdot|_{\can}^d$ on $\mathcal{L}^{d,\an}_{\C((t))}$, which will also be denoted later  by $|\cdot|_{\can}$ through abuse of notation.

Using the canonical metric, we associate the model function to a metric on $\mathcal{L}_{\C((t))}^{\an}$ by
\begin{align*}
    |\cdot|_{\F}:=|\cdot|_{\can}e^{-\frac{1}{d}g_{\F}}.
\end{align*}

In the same way, any function $g:\pr^{2,\an}_{\C((t))}\to\R$ defines a metric on $\mathcal{L}$ via $|\cdot|_g=|\cdot|_{\can}e^{-g}$. In particular, if $g$ is a uniform limit of $\{g_n/d_n\}$ for some sequence of model functions $g_n$ of degree $d_n$, then it also defines a metric, and we can consider its Monge-Amp\`ere measure: for such $g$, define the Monge-Amp\`ere measure $\MA_0(g)$ associated to $g$ is defined by $(\ddc\log|\cdot|_{g})^{\wedge 2}$. The Laplacian operator and the wedge product over Berkovich spaces were developped in \cite{CLD}, and the uniform function $g$ has an interpretation via the language of the intersection theory as in the complex case; for the details, see section 5.6.5 in \cite{CLD}.


\subsection{Basic notions of hybrid spaces}
First, fix a real number $0<r<1$. Next, the hybrid norm $|\cdot|_{\hyb}$ on the complex number field $\C$ is defined by
\begin{align*}
    |x|_{\hyb}=
    \begin{cases}
    \max\{|x|,1\} &\text{ if }x\neq0, \\
    0 &\text{ if }x=0
    \end{cases}
\end{align*}
for any $x\in\C$, where $|x|$ is the Euclidean norm of $x$. Given the hybrid norm on $\C$, define the following ring $A_r$:
\begin{align*}
    A_r:=\left\{f(t)=\sum_{i=-\infty}^{\infty}a_it^i\ \middle|\ a_i\in\C,\ \sum_{i=-\infty}^{\infty}|a_i|_{\hyb}r^i<\infty\right\}.
\end{align*}
We can naturally introduce a metric structure on $A_r$ as follows:
\begin{align*}
    \left\|\sum_{i=-\infty}^{\infty}a_it^i\right\|:=\sum_{i=-\infty}^{\infty}|a_i|_{\hyb}r^i<\infty.
\end{align*}
With this metric, we regard $A_r$ as a completely normed $\C$-algebra, where the norm on $\C$ is the hybrid norm. Then, by taking the Berkovich spectrum of $A_r$, we have an analytic space $\mathcal{C}_{\hyb}(r)$. 

By definition, the points of the Berkovich spectrum of $A_r$ are bounded multiplicative seminorms on $A_r$. Hence, define the map $\tau:\bar{\mathbb{D}}_r\to \mathcal{C}_{\hyb}(r)$ as follows:
\begin{align*}
    \tau(t)(f):=
    \begin{cases}
    |f(t)|^{\frac{\log r}{\log|t|}} &\text{ if }0<|z|\leq r, \\
    |f|_r=r^{\ord_0(f)} &\text{ if }t=0,
    \end{cases}
\end{align*}
for $f\in A_r$ and $t\in\bar{\mathbb{D}}_r$. It is known that $\tau$ is actually a homeomorphism.

The Berkovich analytification of the scheme $\Proj A_r[X,Y,Z]$ is called the hybrid projective plane and is denoted by $\pr^2_{\hyb}$. From \cite{Fa}, we have the following properties:
\begin{itemize}
    \item The natural map $\pi_{\hyb}:\pr^2_{\hyb}\to\mathcal{C}_{\hyb}(r)$ is continuous and proper.
    \item For any $t\neq0$, we have a canonical homeomorphism
    \begin{align*}
        \psi_t:\pr^2(\C)\to\pi_{\hyb}^{-1}\{\tau(t)\}.
    \end{align*}
    \item For $t=0$, there is a canonical homeomorphism
    \begin{align*}
        \psi_0:\pr^{2,\an}_{\C((t))}\to\pi^{-1}\{\tau(0)\}.
    \end{align*}
\end{itemize}
Here, $\pr^{2,\an}_{\C((t))}$ denotes the Berkovich analytification of $\pr^2_{\C((t))}$, where $\C((t))$ is regarded as a complete non-archimedean field by the $t$-adic norm.

Any point $z\in\pr^2_{\hyb}$ defines a bounded norm on a residue field of some (scheme-theoretic) point of $\pr^2_{A_r}$. As $A_r$ is a $\C$-algebra, we can restrict the bounded norm on the residue field to $\C$. Since it is bounded by $|\cdot|_{\hyb}$ by definition, there exists a real number $\eta(z)\in[0,1]$ such that $|\cdot|_z=|\cdot|^{\eta(z)}$ on $\C$, where $|\cdot|_z$ is the seminorm corresponding to $z$ and $|\cdot|$ is the Euclidean norm on $\C$. Note also that $|\cdot|^0$ is the trivial norm on $\C$, and that the function $\eta$ is known to be continuous and proper. 

Also, by taking the Berkovich analytification of $\mathbb{A}^{2}_{A_r}=\Spec A_r[x,y]$, we obtain the hybrid affine plane $\mathbb{A}^2_{\hyb}$, for which properties similar to those listed above for $\pr^2_{\hyb}$ also hold.

\subsection{Monge-Amp\`ere measure on the hybrid projective space}
For any regular admissible datum $\F$, define a measure $\mu_{t,\F}$ for $t\in\bar{\D}_r$ on $\pi^{-1}(t)\subset\pr^2_{\hyb}$ as follows:
\begin{align*}
    \mu_{t,\F}&=(\psi_t)_*(\MA_t(\phi_{\F}))\text{ for }t\neq0, \\
    \mu_{0,\F}&=(\psi_0)_*(\MA_0(g_{\F}))\text{ for }t=0.
\end{align*}

For the family of measures $\{\mu_{t,\F}\}_{t\in\bar{\D}_r}$ on $\pr^2_{\hyb}$, we have a weak convergence
\begin{align*}
    \lim_{t\to0, t\neq0}\mu_{t,\F}=\mu_{0,\F}.
\end{align*}


The most important fact from \cite{Fa} is that the same weak convergence also follows for a uniform function;

\begin{prop}[\cite{Fa}, Definition 4.1 \& Theorem 4.2] \label{prop:themainfact}

Let $\phi:\pr^2(\C)\times\D^*\to\R$ be a uniform function; that is, there exists a sequence of regular admissible data $\{\F_n\}_n$ and a decreasing sequence of positive real numbers $\{\epsilon_n\}_n$ converging to $0$ such that, for any $t$ satisfying $0<|t|\leq r$,
\begin{align*}
     \sup_{z\in \pr^2(\C)\times\{t\}}\left|\phi(z)-\frac{1}{d_n}\phi_{\F_n}(z)\right|\leq \epsilon_n \log|t|^{-1}.
\end{align*}
Then, the sequence of functions $\{g_{\F_n}/d_n\}_n$ on $\pr^{2,\an}_{\C((t))}$ uniformly converges to a continuous function $g$. Moreover, by setting
\begin{align*}
    \MA_{t,\hyb}(\phi)&:=(\psi_t)_*\MA_t(\phi) \text{ for }t\neq0,\\
    \MA_{0,\hyb}(\phi)&:=(\psi_0)_*\MA_0(g) \text{ for }t=0,
\end{align*}
we obtain the following weak convergence:
\begin{align*}
    \lim_{t\to0}\MA_{t,\hyb}(\phi)=\MA_{0,\hyb}(\phi).
\end{align*}
\end{prop}

\section{H\'enon map over the hybrid affine plane} \label{sec:maintheorem}
as in Section \ref{sec: Aoneparameterfamilyofcomplexhenonmaps}, consider an analytic family of H\`enon maps that are parametrized by $\D^*$ and meromorphically degenerate at the origin:
\begin{align*}
    H_t(x,y)=(p_t(x)-a(t)y,x),
\end{align*}
where $p_t(x)=x^d+a_1(t)x^{d-1}+\cdots+a_d(t)\in\mathcal{O}(\mathbb{D})[t^{-1}][x]$, and $a(t)\in\mathcal{O}(\mathbb{D})[t^{-1}]$. This family can be regarded canonically as an automorphism $\A^2_{\C((t))}\to\A^2_{\C((t))}$. Let $H:\A^{2,\an}_{\C((t))}\to\A^{2,\an}_{\C((t))}$ denote its analytification. Also, for a fixed $0<r<1$, we consider the hybrid affine plane. The family $\{H_t\}_{\D^*}$ also canonically induces 
\begin{align*}
    H_{\hyb}:\A^{2}_{\hyb}\to\A^{2}_{\hyb}.
\end{align*}

This section's goal is to study the relation between the family of dynamics defined by $\{H_t\}$ and the non-archimedean dynamics defined by $H$ over the hybrid space.

As in Section \ref{sec: Aoneparameterfamilyofcomplexhenonmaps}, for each $t\neq0$, define
\begin{align*}
    G_{n,t}(x,y):=\frac{1}{d^n}\log\max\left\{\|H_t^{\circ n}(x,y)\|,\|H_t^{\circ -n}(x,y)\|\right\}.
\end{align*}
As was shown there, the sequence $\{G_{n,t}\}_n$ uniformly converges on the whole $\C^2$. Let $G_t$ denote the limit, and define $\mu_t=(\ddc G_t)^{\wedge 2}/4\pi^2$.

On the other hand, similar objects can be defined for the non-archimedean dynamics defined by $H$. For $z\in\mathbb{A}^{2,\an}_{\C((t))}$, define
\begin{align*}
    G_n(x,y):=\frac{1}{d^n}\log\max\left\{\|H^{\circ n}(x,y)\|,\|H^{\circ -n}(x,y)\|\right\}.    
\end{align*}
In Section \ref{sec:non-archimedean-henon-maps}, we showed that the sequence $\{G_n\}_n$ converges locally uniformly to a continuous function, denoted by $G$. Let $\mu$ be a positive measure defined by $(\ddc G)^{\wedge 2}$. 

For a nonzero parameter $t$, we attach the invariant measure $\mu_t=(\ddc{G_t})^{\wedge 2}/4\pi^2$. Also, since all the coefficients of $H_t$ are elements of the non-archimedean field $\C((t))$, we can define the measure $\mu=(\ddc{G})^{\wedge 2}$. We can push-forward these measures to the hybrid affine plane, and we now state the following main theorem.

\begin{thm}\label{mainthm}
In the above setting, as $t\to0$, we have the following weak convergence of measures:
\begin{align*}
    (\psi_t)_*(\mu_t)\to(\psi_0)_*(\mu)
\end{align*}
\end{thm}

The theorem's proof involves application of Proposition \ref{prop:themainfact}, but the key obstacle is that the affine plane is not projective. We can compactify the affine plane in many ways, but none of them admits the extension of $H_t$ as an endomorphism. Hence, we overcome this problem by considering $H^{\circ n}$ and $H^{\circ -n}$ simultaneously.

Let $H^{\circ n}(x,y)=(H_{1,t}^{(n)}(x,y),H_{2,t}^{(n)}(x,y))$, and consider the following sequence of morphisms:
\begin{align*}
    \tilde{H}^{(n)}_t(x,y):=\left(H_{1,t}^{(n)}(x,y),H_{2,t}^{(n)}(x,y),H_{1,t}^{(-n)}(x,y),H_{2,t}^{(-n)}(x,y)\right).
\end{align*}
For each $t\neq0$, $\tilde{H}^{(n)}$ is a morphism from $\C^2$ to $\C^4$. Also, as before, $\tilde{H}^{(n)}_t$ induces the following morphism:
\begin{align*}
    \tilde{H}^{(n)}:\A^{2,\an}_{\C((t))}\to\A^{4,\an}_{\C((t))}.
\end{align*}
For these morphisms, the following fact holds.

\begin{prop}\label{prop:homogenizationofH}
 Let $\tilde{H}^{(n)}_t$ be as above. Then, for any $n$, $\tilde{H}^{(n)}_t$ extends to a morphism $\pr^2(\C)\to\pr^4(\C)$
\end{prop}
\begin{proof}
Consider $\pr^2_{\C}=\Proj{\C[X,Y,Z]}$ such that $x=X/Z$ and $y=Y/Z$. We can naturally extend the morphism $\tilde{H}^{(n)}_t$ by defining
\begin{align*}
    F_t^{(n)}(X:Y:Z)=\left(F_{1,t}^{(n)}(X,Y,Z):F_{2,t}^{(n)}(X,Y,Z):F_{1,t}^{(-n)}(X,Y,Z):F_{2,t}^{(-n)}(X,Y,Z):Z^{d^n}\right),
\end{align*}
where $F_{i,t}^{(n)}$ is the homogenization of $H_{i,t}^{(n)}$. In general, a map $F_t^{(n)}$ is a rational map. To show that it is a morphism, it is sufficient to check that there is no common zero, i.e., that $F_{1,t}^{(n)}(X,Y,Z)=F_{2,t}^{(n)}(X,Y,Z)=F_{1,t}^{(-n)}(X,Y,Z)=F_{2,t}^{(-n)}(X,Y,Z)=Z^{d^n}=0$ only if $X=Y=Z=0$. 

Note that for a natural number $n$, $H_{1,t}^{(n)}$ is always of the form
\begin{align}
    H_{1,t}^{(n)}(x,y)=x^{d^n}+\text{(l.o.t.)}, \label{eq:binomialtheorem}
\end{align}
where the lower-order terms are of the form $c(t)x^iy^j$ with $i+j<d^n$. This can be shown by induction. For $n=1$, it holds by definition. Now, assume the statement is true for any $i\leq n-1$, and let $H_{1,t}^{(i)}=x^{d^{i}}+q_{i,t}(x,y)$. Then,
\begin{align*}
    H^{(n)}_{1,t}(x,y)=\left(x^{d^{n-1}}+q_{n-1,t}(x,y)\right)^d+a_1(t)\left(x^{d^{n-1}}+q_{n-1,t}(x,y)\right)^{d-1}&+\cdots+a_d(t)\\ &-a(t)\left(x^{d^{n-2}}+q_{n-2,t}(x,y)\right),
\end{align*}
from which equation (\ref{eq:binomialtheorem}) follows.

We have a similar statement for $H^{\circ -n}$: for a natural number $n$, the polynomial $H_{2,t}^{(-n)}$ is always of the form
\begin{align*}
    H_{2,t}^{(-n)}(x,y)=y^{d^n}+\text{(l.o.t.)}.
\end{align*}
Hence, taking the homogenization, we have
\begin{align*}
    F_{1,t}^{(n)}(X:Y:Z)=X^{d^n}+ZP^{(n)}_t(X:Y:Z),
\end{align*}
and
\begin{align*}
    F_{2,t}^{(n)}(X:Y:Z)=Y^{d^n}+ZQ^{(n)}_t(X:Y:Z),
\end{align*}
for some homogeneous polynomials $P_t^{(n)}$ and $Q_t^{(n)}$ of degree $d^n-1$. Also, the sum of the degrees of $H_{2,t}^{(n)}$ with respect to $x$ and that with respect to $y$ is strictly smaller than $d^n$, and the same holds for $H_{1,t}^{(-n)}$. Accordingly, we can set $H_{2,t}^{(n)}=ZR_t^{(n)}$ and $H_{2,t}^{(n)}=ZS_t^{(n)}$.  By combining these, we obtain
\begin{align*}
    F^{(n)}(X:Y:Z)=[X^{d^n}+ZP_t^{(n)}:ZR_t^{(n)}:ZS_t^{(n)}:Y^{d^n}+ZQ^{(n)}_t:Z^{d^n]}
\end{align*}
Since the common zero must be in the locus $\{Z=0\}$, we can take $Z=0$ so that we have $F^{(n)}(X:Y:0)=[X^{d^n}:0:0:Y^{d^n}:0]$. Then, the common zero must satisfy $X=Y=Z=0$.
\end{proof}

Next, we regard $F_{1,t}^{(n)}$, $F_{2,t}^{(n)}$, $F_{1,t}^{(-n)}$, and $F_{2,t}^{(-n)}$ as sections of $\mathcal{L}^{\otimes d^n}$. Then, by the above proposition, we have a regular admissible datum $\F_n$ with the following properties:
\begin{itemize}
    \item The degree is $d^n$;
    \item Meromorphic sections are $F_{1,t}^{(n)}$, $F_{2,t}^{(n)}$, $F_{1,t}^{(-n)}$,  $F_{2,t}^{(-n)}$, and $Z^{d^n}$.
    \item $\mathfrak{a}_n$ is generated by $F_{1,t}^{(n)}$, $F_{2,t}^{(n)}$, $F_{1,t}^{(-n)}$, $F_{2,t}^{(-n)}$, and $Z^{d^n}$.
    \item $p_n$ is any log-resolution of $\mathfrak{a}_n$.
\end{itemize}

For convenience in later arguments, define
\begin{align*}
    \Phi_{n,t}(X:&Y:Z):=
    \\&\max\left\{\left|F_{1,t}^{(n)}(X:Y:Z)\right|, \left|F_{2,t}^{(n)}(X:Y:Z)\right|, \left|F_{1,t}^{(-n)}(X:Y:Z)\right|, \left|F_{2,t}^{(-n)}(X:Y:Z)\right|,|Z|^{d^n}\right\}.
\end{align*}

Note that by definition we have
\begin{align*}
    \phi_{\F_n}([X:Y:Z],t)=\log\frac{\Phi_{n,t}([X:Y:Z])}{\left(|X|^2+|Y|^2+|Z|^2\right)^{d^n/2}},
\end{align*}
and
\begin{align*}
    \Phi_n(x:y:1)&=\log\max\left\{\left|F_{1,t}^{(n)}(X:Y:Z)\right|, \left|F_{2,t}^{(n)}(x:y:1)\right|, \left|F_{1,t}^{(-n)}(x:y:1)\right|, \left|F_{2,t}^{(-n)}(x:y:1)\right|,|1|^{d^n}\right\} \\
    &=\log^+\max\left\{\left|H_{1,t}^{(n)}(x,y)\right|,\left|H_{2,t}^{(n)}(x,y)\right|, \left|H_{1,t}^{(-n)}(x,y)\right|, \left|H_{2,t}^{(-n)}(x,y)\right| \right\} \\
    &=G_n(x,y).
\end{align*}
Then, for a sequence of regular models $\{\F_n\}$, we have the folloeing uniformity:

\begin{prop}\label{prop:uniformityofphin}
 Let $\{\F_n\}$ be as above. Then, there exists a continuous function $\phi_{\infty}$ and a decreasing sequence of positive real numbers $\{\epsilon_n\}_n$ converging to $0$ such that
 \begin{align*}
     \sup_{z\in \pr^2(\C)\times\{t\}}\left|\phi(z)-\frac{1}{d^n}\phi_{\F_n}(z)\right|\leq \epsilon_n \log|t|^{-1}
 \end{align*}
 for any $t\in\D^*$.
\end{prop}

\begin{proof}
Note that by definition, for any $(x,y)\in\C^2=\{Z\neq0\}$ with $x=X/Z$ and $y=Y/Z$, we can take a canonical homogenous coordinate $[x:y:1]$ such that
\begin{align*}
    \phi_{\F_n}([x:y:1],t) &=\log\frac{\Phi_{n,t}([x:y:1])}{\left(|x|^2+|y|^2+1\right)^{d^n/2}}\\
    &=\log\max\left\{\left|H_{1,t}^{(n)}(x,y)\right|, \left|H_{2,t}^{(n)}(x,y)\right|, \left|H_{1,t}^{(-n)}(x,y)\right|, \left|H_{2,t}^{(-n)}(x,y)\right|\right\}-\frac{d^n}{2}\log(|x|^2+|y|^2+1),
\end{align*}
i.e., $\phi_{\F_n}([x:y:1],t)/d^n=G_{n,t}(x,y)$. Therefore, the uniformity was already shown in Proposition \ref{prop:mainevaluation}.

It remains to prove the claim for the locus $\{Z=0\}$. However, the proof of Proposition \ref{prop:homogenizationofH} already indicated that
\begin{align*}
    F_t^{(n)}(X:Y:0)=[X^{d^n}:0:0:Y^{d^n}:0].
\end{align*}
Therefore, for any $n$, we have
\begin{align*}
    \phi_{\F_n}([X:Y:0];t)=d^n\log\max\{|X|,|Y|\},
\end{align*}
whose uniformity is trivial.
\end{proof}

Now, we prove Theorem \ref{mainthm}.
\begin{proof}
By Propositionx \ref{prop:themainfact} and \ref{prop:uniformityofphin}, we have
\begin{align*}
    \MA_{t,\hyb}(\phi_{\infty})\to\MA_{t,\hyb}(\phi_{\infty})
\end{align*}
weakly as $0\neq t\to0$, where the function $\phi_{\infty}$ is the one we defined in Proposition \ref{prop:uniformityofphin}. Everything here is restricted to an open subspace $\A^2_{\hyb}=\{Z\neq0\}$ of $\pr^{2,\an}$. It is sufficient to show that
\begin{align*}
    &\MA_{t,\hyb}(\phi_{\infty})|_{\{Z\neq0\}}=\mu_t \text{ for } t\neq0, \text{ and} \\
    &\MA_{0,\hyb}(\phi_{\infty})|_{\{Z\neq0\}}=\mu.
\end{align*}

Since the Monge-Amp\`ere measure on the hybrid space is defined as the push-forward of the Monge-Amp\`ere measures on the complex and non-archimedean spaces, it is sufficient to show the equalities on the complex affine plane for $t\neq0$ and on the non-archimedean affine plane for $t=0$ respectively, for which the arguments are parallel.

Recall that by definition (\ref{eq:complexmetricinducedbyfunctions}), every function induces a metric on $\mathcal{L}$ and the Monge-Amp\`ere measure is defined by its Laplacian. For each $t\neq0$, the function $\phi_{\infty}(t,\cdot)$ induces a metric on $\mathcal{L}$, which is denoted by  $|\cdot|_t$. For any homogeneous polynomial $P$ of degree $1$ and $[X:Y:Z]$ with $P(X:Y:Z)\neq0$, we have
\begin{align} \label{eq:metricoftrivializingsection}
    \log|\tau_P(&X:Y:Z)|_{t}=\lim_{n\to\infty}\frac{1}{d^n}\log\frac{|P(X:Y:Z)|^{d^n}}{\Phi_{n,t}([X:Y:Z])},
\end{align}
where $\tau_{P}$ denotes the section in $\mathcal{L}$ that corresponds to the homogeneous polynomial $P$. Indeed, by the construction of $\phi_{\infty}$ and the definition of $|\cdot|_*$,
\begin{align*} 
    \log|\tau_P(&X:Y:Z)|_{t}=\lim_{n\to\infty}\frac{1}{d^n}\left(d^n\log|\tau_{P}(X:Y:Z)|_{*} -\phi_{\F_n}([X:Y:Z])\right) \nonumber \\
    &=\lim_{n\to\infty}\frac{1}{d^n}\left(\frac{|P(X:Y:Z)|^{d^n}}{\left(|X|^2+|Y|^2+|Z|^2)^{d^n/2}\right)}-\log\frac{\Phi_{n,t}([X:Y:Z])}{\left(|X|^2+|Y|^2+|Z|^2\right)^{d^n/2}}\right) \nonumber \\
    &=\lim_{n\to\infty}\frac{1}{d^n}\log\frac{|P(X:Y:Z)|^{d^n}}{\Phi_{n,t}([X:Y:Z])} ,
\end{align*}

Now, we take $P(X:Y:Z)=Z$; then, on the locus $\{Z\neq0\}$, 
\begin{align*}
    \log|\tau_Z(x,y)|_{t}&=\lim_{n\to\infty}\frac{1}{d^n}\log\frac{|Z|^{d^n}}{\Phi_{n,t}([X:Y:Z])}\\
    &=\lim_{n\to\infty}\frac{1}{d^n}\log\frac{1}{\Phi_{n,t}([x:y:1])} \\
    &=-\lim_{n\to\infty}G_{n,t}(x,y)\\
    &=-G_t(x,y),
\end{align*}
where $x=X/Z$ and $y=Y/Z$. By the definition of $\mu_t$, we thus have $(\psi_t)_*\mu_t=\MA_{t,\hyb}(\phi)$.

On the other hand, the same equation as (\ref{eq:metricoftrivializingsection}) also holds for non-archimedean H\'enon maps. That is, let $|\cdot|_{g_{\infty}}$ be the metric induced from $g_{\infty}(z):=\lim_{n\to\infty}g_{\F_n}(z)/d^n$ (which exists by Proposition \ref{prop:themainfact}). 
Then, for any homogeneous polynomial $P$ of degree $1$ with coefficients in $\C((t))$,
\begin{align*}
    \log|\tau_P(&z)|_{g_{\infty}}=\lim_{n\to\infty}\frac{1}{d^n}\left(d^n\log|\tau_P(z)|_{\can} -g_{\F_n}(z)\right)\\
    &=\lim_{n\to\infty}\frac{1}{d^n}\log\frac{|P(z)|^{d^n}}{\max\left\{\left|F_{1}^{(n)}(z)\right|, \left|F_{2}^{(n)}(z)\right|, \left|F_{1}^{(-n)}(z)\right|, \left|F_{2}^{(-n)}(z)\right|,\left|Z(z)\right|\right\}}.
\end{align*}
Here, the $F_i^{(n)}$ are the same polynomials as $F_{i,t}^{(n)}$ but are regarded as polynomials with coefficients in $\C((t))$, while $\tau_P$ is the section of $\mathcal{O}(1)^{\an}$ that corresponds to $P$, as before. For the definition of $|P(z)|$, see Section \ref{section:non-archimedeanmongeamperemeasure}. Finally, on the locus $\{Z\neq0\}$, as in the complex case, we have
\begin{align*}
    \log|\tau_Z(z)|_{g_{\infty}}=-\lim_{n\to\infty}G_n(z)=-G(z),
\end{align*}
thus yielding $(\psi_0)_*(\mu)=\MA_{0,\hyb}(\phi)$.
\end{proof}
Lastly, we discuss the limit of Lyapunov exponents. As above, we consider an analytic family $\{H_t\}_{t\in\bD^*}$ of H\'enon maps and a family $\{\mu_t\}_t$ of the attached measures. As written in section \ref{sec:intro}, for each $t$, the Jacobi matrix $\Diff H_t$ is written as 
\begin{align} \label{eq:DH}
    \Diff H= \begin{pmatrix}
    p_t'(x) & a(t) \\
    1 & 0
    \end{pmatrix} .
\end{align}
The total Lyapunov exponent $\Lyap(H_t)$ is just $\log|a(t)|$. By the definition of hybrid spaces and the way we regard $\bC^2\times\bD^*$ as a subspace of $\bA^2_{\hyb}$, we trivially have
\begin{align}\label{eq:theconvergenceoftotalLyap}
    \Lyap(H_t)=\frac{\log|t|^{-1}}{\log r}\log|a|+o(\log|t|^{-1}),
\end{align}
where $|a|$ is the norm of $a$ on $\bC((t))$. 

We next consider the first and second Lyapunov exponents. For each $t$, by Oseledets theorem, there exist complex numbers $\lambda_{1,t}$ and $\lambda_{2,t}$ such that $|\lambda_{1,t}|\geq|\lambda_{2,t}|$, and for $\mu_t$-almost everywhere $(x,y)$, the eigenvalues of $DH_t(x,y)$ are $\lambda_{1,t}$ and $\lambda_{2,t}$. Since we have $\Lyap(H_t)=\lambda_{1,t}+\lambda_{2,t}$ by definition, we discuss the first (maximal) Lyapunov exponent $\lambda_{1,t}$ here.

For the induced non-archimedean H\'enon map $H$, set $\mu:=\ddc \max\{G^+,G^-\}^{\wedge 2}$ as above. We have $(\psi_t)_*\mu_t\to (\psi_0)_*\mu$ weakly as $t\to0$ over the hybrid space $\bA^2_{\hyb}$ in Theorem \ref{mainthm}. By means of it, we can show the following theorem:

\begin{thm}
    Let $\{H_t\}_{t\in\bD^*}$ and $\lambda_{1.t}$ be as above. Then,
    \begin{align*}
        \lambda_{1,t}=\frac{\log|t|^{-1}}{\log r}\Lambda+o(\log|t|^{-1}),
    \end{align*}
    where
    \begin{align*}
        \Lambda=\lim_{n\to\infty}\frac{1}{n}\int_{\bA^{2,\an}_{\bC((t))}}\log\|\Diff H^n\|d\mu.
    \end{align*}
\end{thm}

Theorem \ref{thm:convoflyapunocexponents} is combination of the equation (\ref{eq:theconvergenceoftotalLyap}) and this theorem.

\begin{proof}
For $(x,y,t)\in\bC^2\times\bD^*$, let $\chi_t(x,y)(z)\in \bC[z]$ be the characteristic polynomial of $DH_t(x,y)$. By (\ref{eq:DH}), we have
\begin{align*}
    \chi_t(x,y)(z)=z^2-p_t'(x)z+a(t),
\end{align*}
where $p_t'(x)=\partial p_t(x)/\partial x$ is a partial derivative with respect to $x$.
For each $t\in\bD^*$,
\begin{align*}
    X_t:=(\Spec \bC[x,y,z]/(\chi_t))
\end{align*}
is a (perhaps singular) complex variety such that the canonical morphism $\pro_t:X_t\to\bA^2_{\bC}$ is finite. Also, set
\begin{align*}
    X:=\{(x,y,z,t)\in \bC^3\times \bD^*\ |\ \chi_t(x,y)(z)=0\}=\bigcup_{t\in\bD^*}X_t(\bC),
\end{align*}
and $\pro:X\to\bC^2\times\bD^*$. Define a continuous function $\phi:X\to\bR$ by
\begin{align*}
    \phi_t(x,y,z):=\log\max\{|z|,|z-p'_t(x)|\}.
\end{align*}
Note that $\chi_t(x,y)(z)=\chi_t(x,t)(z-p'_t(x))=0$ on whole $X$ by definition. For any $\zeta=(x,y,z,t)$ and $\zeta'=(x,y,z',t)\in X$ such that $\pro(\zeta)=\pro(\zeta')$, $\phi_t(x,y,z)=\phi_t(x,y,z')$ is the maximum of the absolute values of the eigenvalues of $\Diff H$ at $\pro(\zeta)$, i.e., the induced map $\bar{\phi}(x,y,t):\bC^2\to\bR$ gives the first Lyapunov exponent $\lambda_{1,t}$ of $H_t$ for all $t\in\bD^*$ and $\mu_t$-almost $(x,y)\in\bC^2$.

Since the coefficients of $\chi_t(x,y)(z)$ can be regarded as elements of both $A_r$ and $\bC((t))$, we can take the analytification of the base extension to have analytic spaces $X_{\hyb}$ and $X^{\an}_{\bC((t))}$ and analytic morphisms $\pro_{\hyb}:X_{\hyb}\to\bA^2_{\hyb}$ and $\pro_0:X^{\an}_{\bC((t))}\to\bA^{2,\an}_{\bC((t))}$. Define the function $\Phi:X_{\hyb}\to\bR$ by
\begin{align*}
    \Phi(\zeta)=
    \begin{cases}
    \log r\cdot\phi_t(x,y)(z)/\log|t|^{-1}    &  \text{ if } \zeta=(x,y,z,t)\in X,\\
          \log\max\left\{|z(\zeta)|,|(z-p_t'(x))(\zeta)|\right\} & \text{ if } \zeta\in X^{\an}_{\bC((t))}.
    \end{cases}
\end{align*}
By the definition of the topology on $X_{\hyb}$, $\Phi$ is continuous on $X_{\hyb}$. Note that for $\zeta,\zeta'\in X_{\hyb}$ such that $\pro_{\hyb}(\zeta)=\pro_{\hyb}(\zeta')$, we have $\Phi(\zeta)=\Phi(\zeta')$. Thus this induces a map $\bar{\Phi}:\bA^{2}_{\hyb}\to\bR$ and for each $t\in\bD^*$, 
\begin{align*}
    \frac{\log|t|^{-1}}{\log r}\bar{\Phi}(x,y,t)=\bar{\phi}_t(x,t)
\end{align*}
gives $\lambda_{1,t}$ for $\mu_t$-almost $(x,y)\in\bC^2$ via the homeomorphism $\psi_t:\bC^2\to\bA^2_{\hyb}$. If we show that the function $\bar{\Phi}$ is continuous on $\bA^2_{\hyb}$, the weak convergence $(\psi_t)_*\mu_t\to (\psi_0)_*\mu$ implies
\begin{align*}
    \lim_{t\to 0}\frac{\log r}{\log|t|^{-1}}\lambda_{1,t}&=\lim_{t\to 0}\frac{\log r}{\log|t|^{-1}}\int_{\bC^2}\bar{\phi}_t(x,y)d\mu_t\\
    &=\int_{\bA^{2,\an}_{\bC((t))}}\log\max\left\{|z(\zeta)|,|(z-p_t'(x))(\zeta)|\right\}d\mu \\
    &=\lim_{n\to\infty}\frac{1}{n}\int_{\bA^{2,\an}_{\bC((t))}}\log\|DH^n\|d\mu,
\end{align*}
which shows the claim. Hence it remains to show the continuity of $\bar{\Phi}$. Since the map $\pro$ is open on $X$, the continuity on $X$ is straightforward. Also, by Lemma 3.2.4 of \cite{Ber1}, The map $\pro_{\hyb, 0}:X^{\an}_{\bC((t))}\to\bA^{2,\an}_{\bC((t))}$ is open, too, which yields the continuity of $\bar{\Phi}|_{X_{\hyb,0}}$. It is sufficient to show that for any net $\{\xi_n\}_n$ of $\psi(\bC^2\times\bD^*)$ such that $\xi_n\to \xi_0\in \bA^2_{\hyb,0}$ in $\bA^2_{\hyb}$, there exists a net $\{\zeta_n\}_n$ such that $\pro_{\hyb}(\zeta_n)=\xi_n$ for all $n$ and a subnet of it converges to $\zeta_0$ satisfying $\pro_{\hyb}(\zeta_0)=\xi_0.$ As the map $\pro_{\hyb}$ is surjective, we can take $\zeta_n$ such that $\pro_{\hyb}(\zeta_n)=\xi_n$.

First, we show that the set $\{\zeta_n\}_n$ is bounded in $X_{\hyb}\subset\bA^3_{\hyb}$. Define the function $\|\cdot\|:X_{\hyb}\to\bR_{\geq0}$ by
\begin{align*}
    \|\zeta\|:=\max\{|x(\zeta)|,|y(\zeta)|,|z(\zeta)|\}.
\end{align*}
Here, the boundedness of the net $\{\zeta_n\}$ means that the set $\{\|\zeta_n\|\}$ is bounded. Note that as $\{\xi_n\}$ is convergent, $\{\|\xi_n\|\}$ is bounded for sufficiently large $n$. More precisely, we can define the similar function $\|\cdot\|$ on $\bA^2_{\hyb}$ by
\begin{align*}
    \|\xi\|:=\max\{|x(\xi)|,|y(\xi)|\},
\end{align*}
and the continuity of it yields
\begin{align*}
    \|\xi_n\|<\|\xi_0\|+1
\end{align*}
for all $n$ sufficiently large. Set $R:=\|\xi_0\|+1\geq 1$. By the definition of $\pro_{\hyb}$, we have $|x(\xi_n)|=|x(\zeta_n)|$ and $|y(\xi_n)|=|y(\zeta_n)|$, i.e., $|x(\zeta_n)|<R$ and $|y(\zeta_n)|<R$ for any sufficiently large $n$. To show the boundedness of $\{\zeta_n\}$, it is sufficient to show that $\{|z(\zeta_n)|\}$ is bounded for any sufficiently large $n$. Note that for any $\zeta\in X_{\hyb}$,
\begin{align*}
    |(z^2+p'_t(x)z+a(t))(\zeta)|=0.
\end{align*}
Assume $|z(\zeta)|>1$; if not, we have $|z(\zeta)|\leq1<R$. Then the above equation assures that $|z(\zeta)|$ cannot exceed $|p'_t(x)(\zeta)|+|a(t)(\zeta)|$ by the triangle inequality. Since $p'_t(x)$ is of degree $(d-1)$, there exists $M>0$ such that $M$ is independent of the choice of $\zeta$ and $|p'_t(x)(\zeta)|<M|x(\zeta)|^{d-1}$. Since $a(t)\in A_r$, $|a(t)(\zeta)|\leq |a(t)|_{\hyb}$ where $|\cdot|_{\hyb}$ is the norm on $A_r$. Combining them together, we have
\begin{align*}
    |z(\zeta)|\leq M|x(\zeta)|^{d-1}+|a(t)|_{\hyb}<M_1|x(\zeta)|^{d-1}
\end{align*}
for some $M_1>M>0$. Therefore, for $n$ so large that $\|\xi_n\|\leq R$,
\begin{align*}
    \|\zeta_n\|&=\max\{|x(\zeta_n)|,|y(\zeta_n)|,|z(\zeta_n)|\} \\
    &\leq\max\{R, M_1\cdot R^{d-1}\} \\
    &=M_1\cdot R^{d-1}.
\end{align*}
Therefore we have $\|\zeta_n\|<\max\{R, M_1 R^{d-1}\}$, i.e., $\{\zeta_n\}$ is bounded.

Since $\bP^3_{\hyb}$ is compact, so is the set of the form $\{\|\zeta\|<R\}$, i.e. the net $\{\zeta_n\}$ has a convergent subnet. Since $\bA^2_{\hyb}$ is Hausdorff, elementary topological argument shows that the limit $\zeta_0$ of any convergent subnet of $\{\zeta_n\}$ satisfies $\pro(\zeta_0)=\xi_0$, which shows the continuity of $\bar{\Phi}$. Hence the proof is complete.

\end{proof}

\end{document}